\def\year{2022}\relax
\documentclass[letterpaper]{article} % DO NOT CHANGE THIS
\usepackage{aaai22}  % DO NOT CHANGE THIS
\usepackage{times}  % DO NOT CHANGE THIS
\usepackage{helvet}  % DO NOT CHANGE THIS
\usepackage{courier}  % DO NOT CHANGE THIS
\usepackage[hyphens]{url}  % DO NOT CHANGE THIS
\usepackage{graphicx} % DO NOT CHANGE THIS
\usepackage{natbib}  % DO NOT CHANGE THIS AND DO NOT ADD ANY OPTIONS TO IT
\usepackage{caption} % DO NOT CHANGE THIS AND DO NOT ADD ANY OPTIONS TO IT
\DeclareCaptionStyle{ruled}{labelfont=normalfont,labelsep=colon,strut=off} % DO NOT CHANGE THIS
\usepackage{algorithm}
\usepackage{algorithmic}

\usepackage{newfloat}
\usepackage{listings}
\floatstyle{ruled}
\newfloat{listing}{tb}{lst}{}
\floatname{listing}{Listing}

\usepackage{amsmath} 
\usepackage{amsthm}
\usepackage{amssymb}
\usepackage{mathtools}
\usepackage{booktabs}
\usepackage[utf8]{inputenc} % allow utf-8 input
\usepackage{url}            % simple URL typesetting
\usepackage{booktabs}       % professional-quality tables
\usepackage{amsfonts}       % blackboard math symbols
\usepackage{nicefrac}       % compact symbols for 1/2, etc.
\usepackage{microtype}      % microtypography
\usepackage{xcolor}         % colors
\usepackage{multirow}
\usepackage{pdflscape} 
\usepackage{bm}
\usepackage[font=small,labelfont=bf,tableposition=top]{caption}
\DeclareCaptionLabelFormat{andtable}{#1~#2  \&  \tablename~\thetable}
\usepackage{pgfplots}
\usepackage{pgfplotstable}
\usetikzlibrary{matrix}
\usepgfplotslibrary{groupplots}
\pgfplotsset{compat=newest}
\usepackage{subcaption}

\definecolor{marine}{RGB}{0,32,96}
\definecolor{navy}{RGB}{0,0,128}
\definecolor{maroon}{RGB}{128,0,0}
\definecolor{olivegreen}{RGB}{85,107,47}
\definecolor{gray}{RGB}{102,102,102}
\definecolor{green}{RGB}{131,198,210}
\definecolor{blue}{rgb}{0, 0.4470, 0.7410}
\definecolor{skyblue}{rgb}{0.3010, 0.7450, 0.9330}
\definecolor{purple}{rgb}{0.4940, 0.1840, 0.5560}
\definecolor{orange}{rgb}{0.9290, 0.6940, 0.1250}
\definecolor{brown}{RGB}{161,121,124}
\definecolor{deepblue}{rgb}{0.0, 0.0, 1.0}

\usepackage[final,authormarkup=none,commentmarkup=footnote]{changes}
\definechangesauthor[color=maroon]{AE}
\definechangesauthor[color=red]{YS}
\definechangesauthor[color=blue]{yunzhuang}
\definechangesauthor[color=orange]{XL}
\definechangesauthor[color=green]{ACE}

\title{Enhancing Column Generation by a Machine-Learning-Based \\ Pricing Heuristic for Graph Coloring}

\author {
    Yunzhuang Shen\textsuperscript{\rm 1},
    Yuan Sun\textsuperscript{\rm 2},
    Xiaodong Li\textsuperscript{\rm 1},
    Andrew Eberhard\textsuperscript{\rm 3},
    Andreas Ernst\textsuperscript{\rm 4}
}
\affiliations {
    \textsuperscript{\rm 1} School of Computing Technologies, RMIT University, Australia\\
    \textsuperscript{\rm 2} School of Computing and Information Systems, University of Melbourne, Australia\\
    \textsuperscript{\rm 3} School of Science, RMIT University, Australia\\
    \textsuperscript{\rm 4} School of Mathematics, Monash University, Australia\\
    s3640365@student.rmit.edu.au, yuan.sun@unimelb.edu.au,\\
    \{xiaodong.li, andy.eberhard\}@rmit.edu.au,
    andreas.ernst@monash.edu
}

\begin{document}
\maketitle
\begin{abstract}

% Column Generation (CG) is an effective method for solving large-scale Linear Programs (LPs). CG aims to find the columns (i.e., variables) with non-zero values in the optimal solution of an LP, starting with a subproblem containing several columns. Then, CG gradually includes useful columns that can improve the solution of the current subproblem. The useful columns are generated from solving a sequence of pricing problems, which are typically NP-hard. To gain the computational advantage, existing studies use heuristic pricing methods to generate \emph{a few} these columns. Unlike these studies, we propose a Machine Learning based Heuristic Pricing (MLHP) approach that can generate \emph{many} useful columns with affordable computational overhead. Given a pricing problem, our MLHP method leverages an ML model to predict its optimal solution, which is then used to guide a sampling method that can generate an arbitrary number of promising columns. On the graph coloring problem, we empirically show that generating a large number of useful columns by MLHP helps CG better capture columns in the optimal LP solution and it significantly reduces the computational time of CG as compared to other state-of-the-art \added[id=AE]{methods}. The improvement in CG can lead to substantially better performance of B\&P.

Column Generation (CG) is an effective method for solving large-scale optimization problems. CG starts by solving a subproblem with a subset of columns (i.e., variables) and gradually includes new columns that can improve the solution of the current subproblem. The new columns are generated as needed by repeatedly solving a pricing problem, which is often NP-hard and is a bottleneck of the CG approach. To tackle this, we propose a Machine-Learning-based Pricing Heuristic (MLPH) that can \emph{generate many high-quality columns efficiently}. In each iteration of CG, our MLPH leverages an ML model to predict the optimal solution of the pricing problem, which is then used to guide a sampling method to efficiently generate multiple high-quality columns. Using the graph coloring problem, we empirically show that MLPH significantly enhances CG as compared to six state-of-the-art methods, and the improvement in CG can lead to substantially better performance of the branch-and-price exact method.

\end{abstract}

\section{Introduction}

Branch-and-price is a widely-used exact method for solving combinatorial optimization problems~\citep{barnhart1998branch} in the general form of Dantzig–Wolfe reformulation~\citep{vanderbeck2000dantzig}. This formulation often provides a much stronger Linear-Programming relaxation (LP) bound than the more compact formulations of the same problem, which may lead to a significant reduction in the problem's search space. However, solving the LP can be challenging, because it typically has an exponential number of variables (or columns) that cannot be considered all at once.

\begin{figure}
    \centering
	\begin{tikzpicture}[shorten >=1pt,auto,node distance=0.5cm, thick,main node/.style={color = marine, rectangle, fill=gray!10, rounded corners, inner sep=6pt}, scale = 1.0]
	\node[main node] (1) at (0,0) { \footnotesize Restricted Master Problem};
	\node[main node] (2) at (4,0) { \footnotesize Pricing Problem};
	
	\path[every node/.style={color = marine}]
	(1) [color = marine, line width = 0.30mm, ->, bend left = 20] edge node { \footnotesize Dual Solution} (2);
	\path[every node/.style={color = marine}]
	(2) [color = marine, line width = 0.30mm, ->, bend left = 20] edge node { \footnotesize  Columns} (1);

	\end{tikzpicture}
    \caption{Illustration of the iterative process of column generation.}
    \label{fig:cg}
\end{figure}
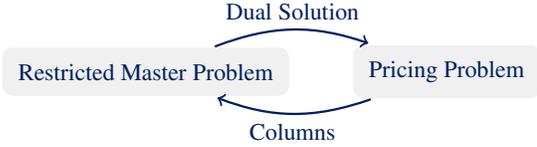

Column Generation~(CG) is an iterative method for solving large-scale LPs. As illustrated in Figure~\ref{fig:cg}, CG starts by solving a subproblem with a small fraction of the columns in an LP, commonly referred to as the Restricted Master Problem (RMP). Then, the optimal dual solution of the RMP is used to set up a pricing problem to search for the column with the least reduced cost. If that column has a negative reduced cost, the column is included in the RMP to further improve its solution. Otherwise, the RMP has captured all the columns with non-zero values in the optimal solution of the original LP. Since an optimal LP solution typically has only a small proportion of columns with non-zero values, CG is expected to solve the LP to optimality without the need to explicitly consider all the columns~\citep{lubbecke2010column}.

In the CG approach, repeatedly solving the pricing problem is typically a bottleneck~\citep{lubbecke2005selected}, because the pricing problem is often NP-hard. To attain computational advantage, heuristic methods are often preferred to trade column quality for computational efficiency. An exact method is usually used only if a heuristic method failed to generate any column with a negative reduced cost~\cite{lubbecke2010column}. Existing studies have explored a variety of pricing heuristics, such as greedy search~\citep{mehrotra1996column, mourgaya2007column} and  metaheuristics~\citep{taillard1999heuristic, malaguti2011exact, beheshti2015novel}. Based on past computational experience, \citet{lubbecke2010column} notes that including multiple columns to the RMP at an iteration of CG can often speed up the progress of CG.

In this paper, we propose a novel Machine-Learning-based Pricing Heuristic (MLPH) for efficiently solving pricing problems. Specifically, we train an ML model offline using a set of solved pricing problems with known optimal solutions. For an unseen pricing problem at an iteration of CG, we use this ML model to predict the optimal solution of the pricing problem, which is then used to guide the search method to generate high-quality columns. To gain efficiency, we employ a linear Support Vector Machine~\citep{boser1992training} for prediction and a sampling method for generating columns. As our method can potentially generate many columns, we introduce several column-selection strategies to form the new RMP to start the next iteration of CG.

By harnessing the knowledge learned from historical data, our MLPH has several advantages over existing pricing methods: (1) compared to sampling-based methods~\citep{dorigo2006ant,cai2016fast}, MLPH can effectively generate columns with \emph{better} reduced costs; (2) compared to other pricing methods~\citep{gurobi2018gurobi, jiang2018two, wang2016two}, MLPH can efficiently generate \emph{many more} columns with negative reduced costs. As our MLPH can \textit{efficiently generate many high-quality columns}, it can help CG capture the columns in an optimal LP solution with a fewer number of iterations.

We demonstrate the efficacy of our proposed MLPH on the graph coloring problem\replaced[id=yunzhuang]{.}{, but we note that our idea can be also applied to other combinatorial optimization problems.} Our experimental study shows that MLPH can significantly accelerate the progress of CG and substantially enhance the branch-and-price exact method.
\section{Background}
In this section, we first introduce different formulations of the Graph Coloring Problem~(GCP). Then, we use GCP to illustrate the solving process of CG. 
\subsection{Graph Coloring Problem}

GCP aims to assign a minimum number of colors to vertices in a graph, such that every pair of the adjacent vertices does not share the same color~\cite{malaguti2010survey}. Let $G(\mathcal{V},\mathcal{E})$ denote a graph, where $\mathcal{V}$ is the set of vertices and $\mathcal{E}$ is the set of edges. GCP can be formulated as: 
 \begin{align}
    \min_{\bm{x},\bm{z}} \;& \sum_{c\in \mathcal{C}}z_c, & \text{(GCP-compact)}\\
        s.t. \;& \sum_{c \in \mathcal{C}} x_{i,c} = 1, &  i \in \mathcal{V},  \\
               & x_{i,c}+x_{j,c} \leq z_c, &  (i,j) \in \mathcal{E}; c \in \mathcal{C},  \\
               & x_{i,c}\in \{0,1\}, \;&  i \in \mathcal{V}; c\in \mathcal{C}, \\
               & z_c \in \{0,1\}, \;&  c\in \mathcal{C}.
            %   & y_{i} \in \mathcal{R};\; i \in \{2, \cdots, n\}
\end{align}
\noindent The binary variable $z_c$ denotes whether a color $c \in \mathcal{C}$ is used to color the graph vertices; and $x_{i,c}$ denotes whether a certain color $c$ is used to color the vertex indexed at $i$. Since this formulation has a polynomial number of variables and constraints, it is commonly referred to as the compact formulation of the GCP, i.e., GCP-compact.

%The objective is to minimize the number of colors used for coloring the graph such that all vertices are colored exactly once (Constraint~2), and the adjacent vertices cannot be associated with the same color (Constraint~3). 

Given that vertices with the same color must be part of an independent set, GCP-compact can be expressed as using a minimum number of Maximal Independent Sets~(MISs) to cover all the vertices in a graph such that every vertex is covered at least once~\citep{mehrotra1996column}, which can be  done systematically using Dantzig–Wolfe decomposition~\citep{vanderbeck2000dantzig,vanderbeck2006generic}. The reformulated problem is commonly referred to as the Set Covering formulation of the GCP~(GCP-SC), defined as:
\begin{align}
    \min_{\bm{x}} \;& \sum_{s \in \mathcal{S}} x_s, & \text{(GCP-SC)}\\
        s.t. \;& \sum_{s \in \mathcal{S},i\in s} x_s \geq 1, &  i \in \mathcal{V},  \\
               & x_s \in \{0,1\}, &s\in \mathcal{S}.
\end{align}
\noindent The binary variable~$x_s$ indicates whether a MIS~$s$ is used to cover a graph, and $\mathcal{S}$ is the set of all the possible MISs in that graph. While GCP-SC provides a much stronger LP than GCP-compact~\citep{mehrotra1996column}, it can contain an exponential number of variables (or columns) to represent all the MISs in a graph. Hence, solving the LP of such a large-scale problem is challenging. 

\subsection{Column Generation}

Given the LP of GCP-SC, CG aims to capture the columns with non-zero values in the optimal LP solution, starting from a RMP with a tiny fraction of the columns in the original LP:
\begin{align}
    \min_{\bm{x_s}} \;& \sum_{s \in \mathcal{\overline{S}}} x_s, & \text{(RMP)}\\
        s.t. \;& \sum_{s \in \mathcal{\overline{S}},i\in s} x_s \geq 1, &  i \in \mathcal{V}, \label{eq:at-least-one} \\
               & 0 \leq x_s \leq 1, &s\in \mathcal{\overline{S}}.
\end{align}
\noindent Note that the integer constraints on $x_s$ are relaxed, and only a small number of MISs is considered initially, i.e., $\overline{\mathcal{S}} \subset \mathcal{S}$. 

\added[id=yunzhuang]{
The RMP can be efficiently solved using the simplex method or the interior point method~\citep{dantzig2016linear}, and its optimal dual solution $\bm{\pi}=[\pi_1,\cdots,\pi_{|V|}]$ associated to vertices (i.e., Constraint~\eqref{eq:at-least-one}) can be used to set up a pricing problem, to search for new MISs with the least reduced cost:
\begin{align}
    \min_{\bm{v}} \;& 1 - \sum_{i \in \mathcal{V}} \pi_i\cdot v_i, & \text{(MWISP)}\\
        s.t. \;& v_i + v_j \leq 1, & (i,j)\in \mathcal{E} \\
              & v_i \in \{0,1\}, & i \in \mathcal{V}.
\end{align}
\noindent The binary variable~$v_i$ denotes whether the vertex~$i$ is a part of the solution, i.e., a MIS according to constraints~(13) and (14). Note that the pricing problem for GCP-SC is the NP-hard Maximum Weight Independent Set Problem (MWISP), where the weight of a vertex $i$ is its dual solution $\pi_i$ to RMP.}

To tackle MWISP, related studies~\cite{mehrotra1996column, malaguti2011exact} employ efficient heuristic methods. Only when a heuristic method fails to find any MIS with a Negative Reduced Cost~(NRC), an exact method is used to solve the MWISP to optimality and so generate the MIS with the least reduced cost. If there exist NRC MISs, they are selectively included in the RMP to further improve its solution, according to a pricing scheme~\citep{lubbecke2005selected}. Otherwise, the RMP has captured all the columns in the optimal solution of the original LP, and hence the original LP is optimally solved.

\section{Machine Learning Based Pricing Heuristic}
\label{sec:MLPH}
Given the MWISP at a CG iteration, we employ an ML model to predict which vertices belong to the optimal MIS. This prediction is then used to guide a sampling method to generate high-quality MISs efficiently. Having many MISs, we introduce several strategies to select a subset of these to form the new RMP at the next CG iteration.

\subsection{Optimal Solution Prediction}
\label{subsec:osp}

We train an ML model to predict the optimal MIS of the MWISP by solving a binary classification task. In our training data, a training example ($\bm{f}$, $y$) corresponds to a vertex in an optimally solved MWISP instance, where~$\bm{f}$ denotes the feature vector that summarizes the property of the corresponding vertex and~$y$ holds a binary value of $1$ (or $0$) indicating whether that vertex is in the optimal~MIS (or not). 

We make use of several features that characterize a vertex of the MWISP, including vertex weight, vertex degree, and the upper bound of a vertex (defined by the sum of weights of that vertex and the vertices that are not adjacent to it). In addition, we adopt two statistical features~\citep{sun2019using} to further enhance the expressiveness of the feature representation for vertices. Given a sample of randomly generated MISs ($\bm{s}\in \mathcal{S}$), the first statistical feature measures the correlation between the presence of a vertex $i$ and the objective values of the sample MISs,
\begin{equation}
\label{eq:corr}
    f_c(i) = \frac{\sum^{K}_{k=1} (\bm{s}^{k}_{i} - \overline{\bm{s}}_{i}) ({o}^{k} - \overline{o})}{\sum^{K}_{k=1} \sqrt{(\bm{s}^{k}_{i} - \overline{\bm{s}}_{i})^2} \sqrt{\sum^{K}_{k=1} ({o}^{k} - \overline{o})^2}},
\end{equation}
where $\bm{s}_i^k$ is a binary value, indicating whether the vertex $i$ is a part of the $k^{th}$ sample; $o^k$ denotes the objective value of that sample; $\overline{\bm{s}}_i$ and $\overline{o}$ respectively denote, the frequency of the vertex $i$ being in a sample and the mean objective value across all samples. A vertex with a high correlation score indicates that this vertex is likely to appear in the high-quality~MISs. 

The second statistical measure uses the rank $r$ of the sample MISs with respect to their objective values,
\begin{equation}
\label{eq:rank}
    f_r(i) = \sum^{K}_{k=1} \frac{\bm{s}^{k}_{i}}{r^{k}}.
\end{equation}
A vertex with a high ranking score indicates that this vertex appears frequently in the high-quality MISs. 

To gain computational efficiency, we adopt Support Vector Machine with linear kernel (linear-SVM)~\citep{boser1992training} to best separate the positive examples (i.e., vertices in the optimal MIS) and negative ones (i.e., vertices not in the optimal MIS) in the training data. For a vertex $i$ in an unseen MWISP instance, the prediction $d_i \in \mathcal{R}$ of the trained linear-SVM is the distance of this vertex, in the feature space, from the optimal decision boundary. This indicates how confidently the linear-SVM classifies this vertex as either in the optimal solution or not according to the signed  distance.

\subsection{Generating Columns via Sampling}
\label{subsec:sampling}

Based on the ML prediction, we can build a probabilistic model to sample multiple high-quality MISs. To generate a MIS, we start with a set containing a randomly selected vertex from the graph, and then iteratively add new vertices into the set until no new vertex can be added. We compute the probability of selecting a vertex via the ML prediction~$d_i$, $i \sim \frac{\sigma(d_i)}{\sum_{j \in \mathcal{C}}^{} \sigma(d_j)};\;i \in \mathcal{C}$, where $\sigma(d_i)$ denotes a logistic function to re-scale the prediction of a vertex into the range of~$[0,1]$, and $\mathcal{C}$ denotes the set of candidate vertices not  already adjacent to any vertex selected. The normalized value~$\sigma(d_i)$ can be interpreted as the `likelihood'  that vertex is in the optimal MIS. Note that we sample MISs starting from a random vertex, so as to increase the diversity of the generated MISs, which has an impact on the solving time of CG.

\subsection{Columns Selection}
\label{subsec:cs}
%  By leveraging the knowledge learned from optimally solved MWISP instances, 
 
 MLPH can potentially generate a large number of NRC columns for unseen MWISPs, and adding all of these to the RMP at early iterations of CG can slow down the solving process of the RMP in the successive CG iterations. However, selectively adding NRC columns may increase the chances of missing out the optimal columns. Therefore, we empirically investigate several strategies to form the RMP at the next CG iteration: 1) \textbf{add-all.} Add all the newly generated NRC columns to the RMP. 2) \textbf{add-partial.} From the newly generated NRC columns, select a proportion of them to add to the RMP in increasing order of their reduced cost. 3) \textbf{replace-existing.} \replaced[id=yunzhuang]{From all the columns, sequentially select columns for the next RMP in the increasing order of their reduced costs, while maintaining the diversity of the set of selected columns by skipping columns too similar to those already added. The algorithm is outlined in the Appendix.}{Replace some of the columns in the current RMP with some of the new columns. The columns with non-zero values in the current RMP solution remain to ensure the solution of the new RMP is no worse than the current. The rest of the columns are examined and included in the increasing order of their reduced cost. To maintain the diversity of the columns, a column is excluded if it was too similar to the columns already included in the new RMP.} 
    
    % \item \textbf{sample-new.} From newly generated columns for the current pricing problem, sample $n$ out of all these columns according to the sampling distribution defined by $p_i=\frac{-c_i}{\sum_{i} -c_i}$, where $p_i$ is the probability of adding the $i^{th}$ column to the RMP among new columns with negative reduced costs $c_i$.
        
    % Note that the same operation is also carried out for the next column selection strategy, i.e., random-all.
    
    % Algorithmically, we create an array, where each entry associates to a graph vertex, and allows to store upto~$k$ columns (i.e., a MIS) containing that vertex, where~$k$ is the ratio between the number of total columns and the number of vertices in a graph. This greedy algorithm goes through the set of candidate columns sorted by their reduced costs in ascending order, and attempt to add this column into an entry (the MIS must contain the vertex corresponding to that entry) with less than~$k$ MISs. If all the entries are full, the greedy algorithm will skip that MIS and continue.

At one end of the spectrum, the strategy~(1) adds all the columns to the RMP, resulting in the fastest growth of the size of the RMP. On the other hand, the strategy~(3) replaces some of the columns in the current RMP with newly generated columns, and it can maintain a fixed number of columns in the RMP. Due to this restriction, CG may require more iterations to capture all the optimal LP columns. The strategy~(2) is somewhere in between these two extremes, resulting in relatively slow growth of the size of the RMP. 
\section{Experiment Settings}

\paragraph{Graph benchmarks and problem instance generation.} We use standard Graph Coloring Benchmarks\footnote{\url{https://sites.google.com/site/graphcoloring/files}}. Given a graph with~$n$ vertices, the goal is to find the columns with non-zero values in the optimal LP solution of GCP-SC, starting from a RMP initialized with $10n$ randomly generated columns. We note that reducing the number of samples can affect the performance of CG negatively. Among $136$ benchmark graphs, we remove those whose initial RMPs already contain all the optimal columns and whose initial RMPs cannot be solved by an LP solver within a reasonable time. For the remaining~$89$ graphs, we label $81$ of them as `small' and $8$ of them as `large', according to the computational time for solving their initial RMPs. For a graph, we can generate multiple RMPs by seeding the initial set of random columns, and these RMPs can be viewed as individual problem instances because solving them can result in different optimal dual solutions and hence different subsequent MWISPs. For training, we generate $10$ instances on $10$ small graphs with random seed $s=1314$. For testing, we generate $24$ instances on each graph using random seeds $s \in \{1,2,\cdots,24\}$, resulting in a total number of $1944$ small instances and $192$ large instances. 

\paragraph{Data collection and training.} For each training instance, we run CG using an exact, specialized solver TSM~\citep{jiang2018two} to solve MWISPs to optimality. The MWISPs with optimal solutions are recorded every five CG iterations up to the $25^{th}$ iteration of CG. In the training data, the statistical features are computed from a set of $n$ MISs, randomly sampled uniformly, with all features normalized instance-wise~\citep{khalil2016learning}. The parameters for training SVM are set to the default values  of~\citep{chang2011libsvm}, except that the regularization term for misclassifying positive training examples is raised to the ratio between the negative training examples and the positive ones. For tuning the parameters in the logistic function, we employ Bayesian Optimization (BO)~\citep{snoek2012practical, nogueira2014bayesian}. Specifically, BO treats the MLPH as a black-box and attempts $300$ runs of MLPH using different sets of parameters in the logistic function to minimize the reduced cost of the best-found solution of the MWISP at the first CG iteration. 

\paragraph{Compared pricing methods.} 
\begin{itemize}
    \item \textbf{Gurobi,} a state-of-the-art commercial Mixed-Integer-Programming (MIP) solver~\citep{gurobi2018gurobi}. Such a MIP solver is used as the pricing method for GCP by related work~\citep{malaguti2011exact}. By default, Gurobi aims to solve a MWISP to optimality, hence it spends most of its computational time on improving the duality bound. In addition to this default configuration, we include another setting, Gurobi-heur, that focuses on finding feasible solutions. This is done by setting the parameters `PoolSearchMode' to $2$, `PoolSolutions' to $10^8$, and `Heuristics' to $95\%$.
    
    \item \textbf{Ant Colony Optimization (ACO),} an efficient meta-heuristic that has been investigated for many combinatorial optimization problems~\citep{dorigo2006ant}. ACO maintains a probabilistic distribution during the solving process and constructs solutions by sampling from that distribution. We adopt the ACO variant as described in \citep{xu2007improved}.
    
    \item \textbf{Specialized methods.}  Since the optimal solution of a MWISP is the same as that of solving the Maximum Weight Clique Problem (MWCP) in its complementary graph, we also include three state-of-the-art MWCP solvers: 1) TSM~\citep{jiang2018two}, an exact solver based on the branch-and-bound framework with domain-specific knowledge for tightening the dual bounds; 2) LSCC~\citep{wang2016two}, a heuristic method based on Local Search; 3) Fastwclq~\citep{cai2016fast}, a heuristic method that constructs solutions in a greedy fashion with respect to a benefit-estimation function. 
\end{itemize}

\paragraph{Computational budgets, evaluation criteria, and other specifications.} Table~\ref{tab:budget} shows the computational budgets for solving small and large problem instances, respectively. In addition to an overall cutoff time for CG, we also set a cutoff time for solving the MWISP at every CG iteration. In particular, the exact methods are also subject to the time limit and are evaluated as heuristic methods. Moreover, we set for each of the pricing methods an individual termination condition. For LSCC based on Local Search (LS), we terminate LSCC if it cannot find better solutions for $50n$ LS iterations. This is because, empirically, LSCC can find high-quality solutions efficiently, and providing excessive computational time can hardly improve these solutions. For MLPH, ACO, and Fastwclq, we set the number of constructed solutions to $50n$. For TSM and Gurobi, they may terminate early when solving a MWISP instance to optimality. We use `add-partial' as the default column-selection method with the column limit $n$. When a pricing method fails to find any NRC column, TSM is used to solve the MWISP to optimality and the optimal column is added to the RMP to start the next iteration. 

\begin{table}[t!]
    \centering
    \resizebox{0.9\columnwidth}{!}{\begin{tabular}{ccccc}
        \toprule
         Label & Total \# Instances & \begin{tabular}{@{}c@{}}Cutoff Time \\ (Overall)\end{tabular} & \begin{tabular}{@{}c@{}}Cutoff Time \\ (Pricing)\end{tabular} & \begin{tabular}{@{}c@{}}\# CPUs  \\ (Paralleled)\end{tabular}\\
         \midrule
         small & 1944 & 1800s & 30s & 1\\
         large & 192 &  8000s & 150s & 4\\
         \bottomrule
    \end{tabular}}
    \caption{Test instances and Computational Budget.}
    \label{tab:budget}
\end{table}

We evaluate the performance of CG using a certain pricing heuristic based on two criteria, the computational time of CG for solved problem instances and the objective value of the RMP (the lower the better) for unsolved problem instances. The latter measures how close the solution of the RMP is to the optimal LP solution and so reflects the progress that CG has made. When reporting the results, we will address CG using a certain pricing method in short, e.g., CG-MLPH.  

During CG the RMPs are solved by the default LP solver of Gurobi~\citep{gurobi2018gurobi}. The experiment is conducted on a cluster with $8$ nodes. Each node has $32$ CPUs (AMD EPYC Processor, 2245 MHz) and $128$ GB RAM. \added[id=yunzhuang]{Our code is written in C/C++ and is available online}\footnote{\url{https://github.com/Joey-Shen/MLPH.git}}.

\section{Results \& Analysis}
\label{exp:cg}

% We first report the results of CG using different pricing methods for solving small problem instances. Then, we report the results for large problem instances, and investigate why a pricing method can help CG converge faster than another. Lastly, we examine the effect of several column selection~strategies.

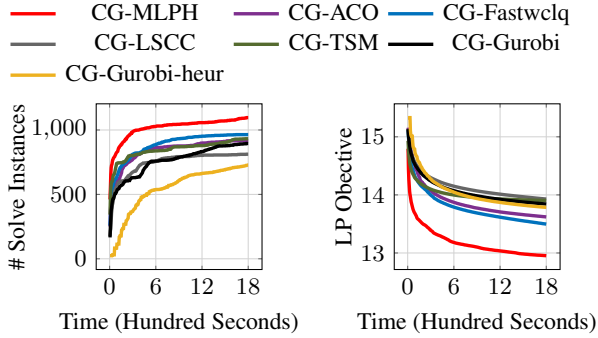
\begin{figure}[t!]
	\begin{tikzpicture}
	\begin{groupplot}[group style = {group size = 2 by 1, horizontal sep = 50pt}, height=0.45\columnwidth, width=0.45\columnwidth, grid style={line width=.1pt, draw=gray!10},major grid style={line width=.2pt,draw=gray!30}, xlabel = \small Time (Hundred Seconds), xtick = {0,6,12,18}, ticklabel style = {font=\small}, xmajorgrids=true, ymajorgrids=true,  major tick length=0.05cm, minor tick length=0.0cm, legend style={font=\small, column sep = 1pt, legend columns = 3,draw=none}]
   \nextgroupplot[%
    legend to name=group1,
    ylabel = \small \# Solve Instances,
    y label style={at={(axis description cs:-0.375,.5)},anchor=south},
    ]

	\addplot[color=red, line width=0.45mm] table [x=x, y=y, col sep=comma] {data/cg/small/solving-curve/svm-cs-0.txt};\addlegendentry{\small CG-MLPH}
	\addplot[color=purple, line width=0.45mm] table [x=x, y=y, col sep=comma] {data/cg/small/solving-curve/aco.txt};\addlegendentry{\small CG-ACO}
	\addplot[color=blue, line width=0.45mm] table [x=x, y=y, col sep=comma] {data/cg/small/solving-curve/fastwclq.txt};\addlegendentry{\small CG-Fastwclq}
	\addplot[color=gray, line width=0.45mm] table [x=x, y=y, col sep=comma] {data/cg/small/solving-curve/lscc.txt};\addlegendentry{\small CG-LSCC}
	\addplot[color=olivegreen, line width=0.45mm] table [x=x, y=y, col sep=comma] {data/cg/small/solving-curve/tsm.txt};\addlegendentry{\small CG-TSM}
	\addplot[color=black, line width=0.45mm] table [x=x, y=y, col sep=comma] {data/cg/small/solving-curve/gurobi.txt};\addlegendentry{\small CG-Gurobi}
	\addplot[color=orange, line width=0.45mm] table [x=x, y=y, col sep=comma] {data/cg/small/solving-curve/gurobi-heur.txt};\addlegendentry{\small CG-Gurobi-heur}
   \nextgroupplot[%
    ylabel = \small LP Obective,
    y label style={at={(axis description cs:-0.2,.5)},anchor=south},
    ]
	\addplot[color=red, line width=0.45mm] table [x=x, y=y, col sep=comma] {data/cg/small/lp-curve/svm-cs-0.txt};
	\addplot[color=purple, line width=0.45mm] table [x=x, y=y, col sep=comma] {data/cg/small/lp-curve/aco.txt};
	\addplot[color=blue, line width=0.45mm] table [x=x, y=y, col sep=comma] {data/cg/small/lp-curve/fastwclq.txt};
	\addplot[color=gray, line width=0.45mm] table [x=x, y=y, col sep=comma] {data/cg/small/lp-curve/lscc.txt};
	\addplot[color=olivegreen, line width=0.45mm] table [x=x, y=y, col sep=comma] {data/cg/small/lp-curve/tsm.txt};
	\addplot[color=black, line width=0.45mm] table [x=x, y=y, col sep=comma] {data/cg/small/lp-curve/gurobi.txt};
	\addplot[color=orange, line width=0.45mm] table [x=x, y=y, col sep=comma] {data/cg/small/lp-curve/gurobi-heur.txt};
    \end{groupplot} 
    \node at (group c1r1.north) [anchor=north, yshift=1.6cm, xshift=1.5cm] {\pgfplotslegendfromname{group1}}; 
	\end{tikzpicture}
\caption{Results for CG with different pricing methods for solving small problem instances. \textbf{Left:} the number of solved instances. \textbf{Right:} the objective values of the RMP (the lower the better), averaged over all problem instances using the geometric mean. }
\label{fig:ret_small}
\end{figure}

% \begin{table}[t!]
%     \caption{Ranking the performances of CG using different pricing methods for small graphs.}
%     \label{tab:small}
%     \centering
%     \resizebox{\columnwidth}{!}{\begin{tabular}{cccccccc}\toprule
%       Methods & MLPH & TSM &  ACO & Fastwclq & LSCC & Gurobi & Gurobi-heur \\ \midrule
%       $1^{st}$ Place & 52 & 12 & 8 & 7 & 1 & 1 & 0 \\
%       $2^{nd}$ Place \& Above & 69 & 26 & 23 & 20 & 1 & 10 & 13 \\
%       $3^{rd}$ Place \& Above & 75 & 39 & 34 & 41 & 5 & 30 & 19 \\
%     %   $4^{th}$ Place \& Above & 78 & 46 & 49 & 54 & 27 & 42 & 28 \\
%         \bottomrule
%     \end{tabular}}
% \end{table}

% Next, we rank the performance of each method with respect to a test graph, and report in Table~\ref{tab:small} the number of times that a method achieves a certain ranking and above. We use the following rule for ranking the performances of CG using different pricing methods: for a graph, a method with lower mean LP objectives over the $24$ problem instances generated using this graph has a better rank. If multiple methods solve all problem instances to optimality (i.e., they have the same LP objective), the method with a shorter average solving time is ranked higher. Table~\ref{tab:small} shows the results. CG using MLPH wins first place on $52$ out of $81$ graphs and it is in the top three pricing methods for $78$ graphs; LSCC and Gurobi-heur are not competitive overall; the rest of the methods are competitive and can be suitable for certain.

\begin{figure}[t!]
    \centering
    \begin{tikzpicture}
    	\begin{axis}[%
    	height=0.5\columnwidth, width=0.9\columnwidth, grid style={line width=.1pt, draw=gray!10},major grid style={line width=.2pt,draw=gray!30}, xlabel = \small Number of Nodes, xtick = {0,100, 1000,2500}, xmajorgrids=true, ylabel= \small Graph Density, ymajorgrids=true, xmode=log, log ticks with fixed point, major tick length=0.05cm, minor tick length=0.0cm, legend style={at={(1.15,1.4)}, column sep = 1pt, legend columns = 4,draw=none}, ticklabel style = {font=\small}, scatter/classes={%
    		MLPH={mark=*,red},
    		ACO={mark=*,purple},
    		Fastwclq={mark=*,blue},
    		LSCC={mark=*,gray},
    		TSM={mark=*,olivegreen},
    		Gurobi={mark=*,black},
    		Gurobi_heur={mark=*,orange}}]
    	\addplot[scatter,only marks,%
    	    mark size=1.5pt,
    		scatter src=explicit symbolic]%
    	table[x=x, y=y,meta=label,col sep=comma] {data/cg/graph_compare.txt};
            \legend{\small CG-MLPH, \small CG-ACO, \small CG-Fastwclp, \small CG-LSCC, \small CG-TSM, \small CG-Gurobi, \small CG-Gurobi-heur}
    	\end{axis}
    \end{tikzpicture} \caption{$81$ small graphs labeled by CG with the winning pricing method. MLPH is the best pricing method for $52$ graphs, followed by TSM for $12$ graphs, ACO for $8$ graphs, Fastwclq for $7$ graphs, Gurobi for $1$ graph, and LCSS for $1$ graph.}
    \label{fig:gcompare}
\end{figure}
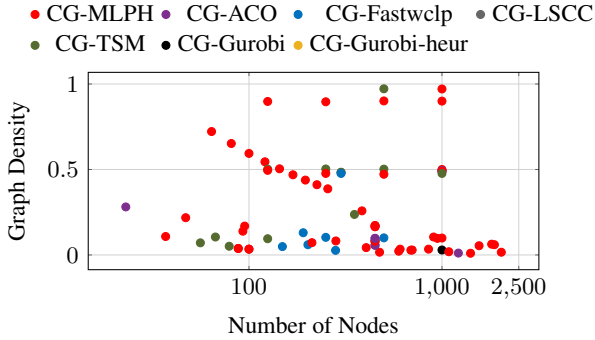

\paragraph{Results for CG using different pricing methods.} Figure~\ref{fig:ret_small} shows the solving statistics for small graphs. The left sub-figure shows that CG-MLPH can solve many more problem instances than CG using other pricing methods (with a given computational budget). The right sub-figure shows that CG-MLPH can make more substantial progress than the comparison methods over all test instances. Noticeably, CG-Fastwclq and CG-ACO have comparable performances, and they are better than the remaining methods we consider. To better understand the CG using different pricing methods with respect to graph characteristics, we report in Figure~\ref{fig:gcompare} the best method on individual graphs. Overall, we observe that MLPH is particularly suitable for CG on relatively larger and/or denser graphs. In contrast, Fastwclp, ACO, and TSM are only competitive on small and sparse graphs. 

Table~\ref{tab:large} shows the results on large graphs for CG. \replaced[id=yunzhuang]{Note that the cutoff time is extended to $8,000$ seconds, and the compared pricing methods are parallelized on $4$ CPUs.}{using several pricing methods that can be paralleled.} It can be seen that CG-MLPH achieves the best performance on every individual graph, and it outperforms the other methods substantially in most cases. In contrast, the other methods are only competitive for certain graphs. Notably, CG-MLPH can optimally solve some problem instances, while other methods cannot solve any instance to optimality. 
\begin{table}[t!]
    \centering
    \resizebox{\columnwidth}{!}{\begin{tabular}{@{}lrlrrrrr@{}}
        \toprule
        Graph & \# Nodes & Density & CG-MLPH & CG-ACO & CG-Gurobi & CG-Gurobi-heur & CG-Fastwclq  \\ 
        \cmidrule(lr){1-3} \cmidrule(lr){4-8}
        wap04a & 5231 & 0.022 & \textbf{44.45} & 75.07 & 56.0 & 55.46 & 81.13 \\
        wap03a & 4730 & 0.026 & \textbf{44.87} & 70.09 & 54.09 & 54.50 & 76.04 \\
        4-FullIns\_5 & 4146 & 0.009 & \textbf{6.68} & 12.32 & 9.60 & 7.67 & 7.91 \\
        C4000.5 & 4000 & 0.500 & \textbf{263.47} & 264.55 & 304.80 & 304.89 & 286.65 \\
        wap02a & 2464 & 0.037 & \textbf{40.02} (3) & 58.57 & 42.69 & 43.49 & 64.52 \\
        wap01a  & 2368 & 0.040 & \textbf{41.0} (24) & 56.21 & 42.71 & 43.34 & 62.13 \\
        C2000.5  & 2000 & 0.500 & \textbf{137.91} & 140.78 & 164.60 & 164.62 & 141.98 \\\
        ash958GPIA & 1916 & 0.007 & \textbf{3.37} & \textbf{3.37} & 3.45 & 3.42 & 3.41 \\
        \cmidrule(lr){1-3} \cmidrule(lr){4-8}
        Geometric Mean & - & - & \textbf{37.09} & 49.20 & 43.08 & 42.15 & 49.49 \\
        \bottomrule
    \end{tabular}}
    \caption{The mean LP objective values (the lower the better) for CG with different pricing methods for large graphs. The results are averaged over the $24$ problem instances generated using the same graph. The number of solved instances (if any) is shown in brackets.}
    \label{tab:large}
\end{table}

\begin{table*}[t!]
    \centering
    \resizebox{0.95\textwidth}{!}{\begin{tabular}{@{}lrlrrrrrrrrrr@{}}
        \toprule
                \multirow{2}{*}{Graph} & \multirow{2}{*}{\# Nodes} & \multirow{2}{*}{Density} & \multicolumn{5}{c}{\# Columns with Negative Reduced Costs} & \multicolumn{5}{c}{Minimum Reduced Cost} \\
         & & & MLPH & ACO & Gurobi & Gurobi-heur & Fastwclq & MLPH & ACO & Gurobi & Gurobi-heur &Fastwclq \\ 
         \cmidrule(lr){1-3}\cmidrule(lr){4-8}\cmidrule(lr){9-13}
        wap04a & 5231 & 0.022 & 191494.5 & 277.8 & 5.5 & 23.0 & 0.0 & -2.48 & -0.37 & -3.11 & -2.94 & N/A \\
        wap03a & 4730 & 0.026 & 234610.1 & 277.8 & 5.8 & 42.1 & 0.0 & -2.44 & -0.35 & -3.05 & -3.04 & N/A \\
        4-FullIns\_5 & 4146 & 0.009 & 28273.2 & 26.3 & 3.0 & 4820.1 & 676.8 & -4.37 & -0.55 & -4.37 & -4.38 & -3.93 \\
        C4000.5 & 4000 & 0.500& 185611.9 & 140979.4 & 2.4 & 1.1 & 89.8 & -0.49 & -0.39 & -0.24 & -0.15 & -0.30 \\
        wap02a & 2464 & 0.037 & 123309.2 & 202.7 & 10.8 & 1064.8 & 0.0 & -1.76 & -0.33 & -2.30 & -2.29 & N/A \\
        wap01a & 2368 & 0.040 & 118508.7 & 243.0 & 10.6 & 1132.1 & 0.0 & -1.79 & -0.36 & -2.32 & -2.33 & N/A\\
        C2000.5 & 2000 & 0.500& 89512.1 & 91193.0 & 1.9 & 1.7 & 253.5 & -0.50 & -0.42 & -0.23 & -0.23 & -0.38\\
        ash958GPIA & 1916 & 0.007 & 95888.6 & 1962.9 & 29.1 & 1507.0 & 58.0 & -0.12 & -0.05 & -0.20 & -0.20 & -0.06 \\
        \bottomrule
    \end{tabular}}
    \caption{Results for different pricing methods solving the MWISP at the initial CG iteration for large problem instances. The first statistic shows the number of columns with negative reduced costs that a pricing method can find. The second statistic shows the reduced cost of the best column that a pricing method can find. Both statistics are averaged over $24$ problem instances generated using the same graph.}
    \label{tab:pp}
\end{table*}

\paragraph{MLPH as a competitive pricing method.} For the set of $8$ large graphs, Table~\ref{tab:pp} compares the performances of different pricing methods for solving the pricing problems in the initial CG iteration (to ensure the results used for compared methods are from solving the same set of MWISPs). Comparing MLPH with other pricing methods, MLPH can find many more NRC columns. Furthermore, the quality of the best-found column by MLPH is highly competitive. Subsequently, CG-MLPH achieves the best performance overall~(Table~\ref{tab:large}). In addition, it can be noted that CG-MLPH solves all problem instances on the graph `wap01a' to optimality using an average number of $36.5$ iterations, while CG methods based on other pricing methods cannot optimally solve any of these, having at a minimum,  an average number of $43.7$ iterations. For other pricing methods, we can also observe that the good performance of CG with a certain pricing method is often accompanied by finding many high-quality columns by that pricing method, such as ACO for dense graphs (e.g., `C4000.5') and Fastwclq for sparse graphs (e.g., 4-FullIns\_5). Similar observations can be also made from the results on small graphs as shown in the Appendix. 

The results indicate that finding a large number of high-quality columns can accelerate the progress of CG. In particular, our proposed MLPH can find a large number of high-quality NRC columns, thereby helping CG obtain much better LP objective values for unsolved problem instances or spending many fewer CG iterations for solved problem instances.

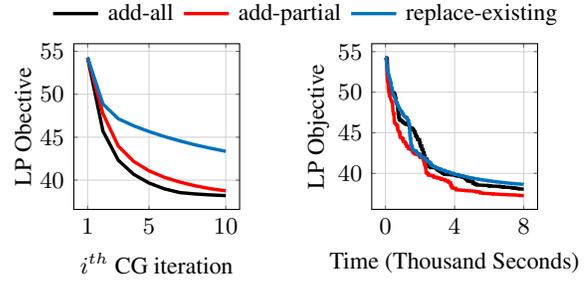
\begin{figure}[t!]
    \centering
	\begin{tikzpicture}
	\begin{groupplot}[group style = {group size = 2 by 1, horizontal sep = 50pt}, height=0.45\columnwidth, width=0.45\columnwidth, grid style={line width=.1pt, draw=gray!10},major grid style={line width=.2pt,draw=gray!30}, xmajorgrids=true, ymajorgrids=true,  major tick length=0.05cm, minor tick length=0.0cm, legend style={font=\small, column sep = 1pt, legend columns = 3,draw=none},ticklabel style = {font=\small}]
   \nextgroupplot[%
    ylabel = \small LP Obective,
        y label style={at={(axis description cs:-0.2,.5)},anchor=south},
        xlabel = \small $i^{th}$ CG iteration, xtick = {1,5,10}
    ]
	\addplot[color=black, line width=0.45mm] table [x=x, y=y, col sep=comma] {data/cg/large-cs/lp-curve-cg/svm-cs-1.txt};
	\addplot[color=red, line width=0.45mm] table [x=x, y=y, col sep=comma] {data/cg/large-cs/lp-curve-cg/svm-cs-0.txt};
	\addplot[color=blue, line width=0.45mm] table [x=x, y=y, col sep=comma]
	{data/cg/large-cs/lp-curve-cg/svm-cs-5.txt};
    \nextgroupplot[%
    legend to name=group2,
    ylabel = \small LP Objective,
    y label style={at={(axis description cs:-0.2,.5)},anchor=south},
    xlabel = \small Time (Thousand Seconds), xtick = {0,4, 8},
    ]
	\addplot[color=black, line width=0.45mm] table [x=x, y=y, col sep=comma] {data/cg/large-cs/lp-curve/svm-cs-1.txt};\addlegendentry{add-all}
	\addplot[color=red, line width=0.45mm] table [x=x, y=y, col sep=comma] {data/cg/large-cs/lp-curve/svm-cs-0.txt};\addlegendentry{add-partial}
	\addplot[color=blue, line width=0.45mm] table [x=x, y=y, col sep=comma] {data/cg/large-cs/lp-curve/svm-cs-5.txt};\addlegendentry{replace-existing}
    \end{groupplot} 
    
    \node at (group c1r1.north) [anchor=north, yshift=0.8cm, xshift=2cm]
    {\pgfplotslegendfromname{group2}}; 
	\end{tikzpicture}
\caption{CG-MLPH using different column-selection strategies for large problem instances (on the left, the $x$-axis is the number of CG iterations; on the right, the $x$-axis is the wall clock time).}
\label{fig:ret_cs_large}
\end{figure}

\paragraph{Efficiency and effectiveness trade-off in column selection.} As shown in the left of Figure~\ref{fig:ret_cs_large}, adding all the NRC columns generated at every CG iteration (i.e., `add-all') results in the faster convergence of the LP objective. On the other hand  keeping a fixed number of columns in the RMP, by replacing the existing columns already in the RMP with newly generated NRC columns (i.e., `replace-existing'), tends to slow down the progress of CG. This shows that adding more columns can increase the chance of capturing the optimal LP columns. When measuring the progress of CG in wall-clock time as shown on the right, we observe that `add-all' and `replace-existing' are comparable because `add-all' increases the computational burden for solving the fast-growing RMP. Compared to these two methods, adding a proportion ($n$ in our case) of the best NRC columns (`add-partial') better balances the trade-off between efficiency and effectiveness.

\section{Branch-and-price with MLPH}

In this part, we use CG-MLPH to enhance Branch-and-Price~(B\&P), an exact method that solves a GCP to optimality by recursively decomposing the original problem (root node) into subproblems (child nodes). During the solving process, CG is used at every node to compute their LP bounds (lower bounds), and a node can be safely pruned without further expansion if its lower bound is no better than the current best-found solution. 

\subsection{Setup}

We use the B\&P code from an open-source MIP solver, SCIP~\citep{gamrath2020scip}. For CG, in particular, an efficient greedy search is used as the pricing heuristic for tackling MWISPs. Only when it fails to generate any NRC column, an exact method called $t$-clique is used to find the optimal column to either certify the optimality of the LP or start the next CG iteration with the optimal column added to RMP. Once the LP at the current node is solved, the node is branched into child nodes, and the columns generated at this node are passed into the child nodes. Furthermore, the B\&P implementation incorporates specialized techniques to GCP from previous studies~\citep{mehrotra1996column,malaguti2011exact}, and more details can be found in the Appendix. 

% In the latter case, the optimal solution to the MWISP is used to compute the Lagrangian lower bound, which is then used to compare with the objective value of the current RMP to detect early termination of CG with optimality guarantee~\cite{malaguti2011exact}

% In the latter case, the optimal solution for the MWISP is then used to compute the Lagrangian lower bound~$L$~\cite{malaguti2011exact}. If $\lceil L \rceil = \lceil O\rceil$ (with $O$ being the objective value of the current RMP), CG is terminated earlier and $\lceil L \rceil$ is set as a valid lower bound of the current node. Then, the current node is branched into child nodes, with the columns generated at this node also passed into the child nodes.  

% the root node, at most $n$ 

We refer to the default setup of B\&P as B\&P-def, and compare it with B\&P-MLPH that replaces the greedy search in B\&P-def with the MLPH for solving MWISPs. Although the greedy search has negligible computational cost, it is less effective as it can only construct a single column at a CG iteration. In contrast, our MLPH can sample many high-quality columns effectively. Empirically, we examine the sample size $\lambda$ of MLPH in $\{10n, n, 0.1n\}$, and observe that no single sample size can fit all graph benchmarks. To best contrast our B\&P-MLPH with B\&P-def, we report the results when $\lambda=10n$. From newly generated NRC columns, we add at most $\theta$ columns into the RMP in the increasing order of their reduced costs, $\theta=n$ for the root node and $\theta=0.1n$ for child nodes. We observe that setting a smaller $\theta$ in column selection for child nodes can reduce the memory required for storing all the columns generated during the B\&P process without sacrificing the performance. This is because child nodes often have a sufficient number of quality columns inherited from their parents and their initial RMPs are already close to the optimum.

For each method, a total number of $1584$ seeded runs are performed to solve the GCPs on a set of~$66$ graphs in the Graph Coloring Benchmarks. The excluded graphs are either too easy (both methods can solve them within $10$ seconds) or too hard to solve (both methods cannot solve the LP at the root node) under the cutoff time of $8000$ seconds. When the LP at the root node is solved to optimality, we report optimality gap, defined as $Gap = 100\% \times \frac{upper\_bound - global\_lower\_bound}{upper\_bound}
$. The upper bound is the objective value of the best-found solution and the global lower bound is determined by the smallest lower bound amongst the remaining open tree nodes. 

\subsection{Results}

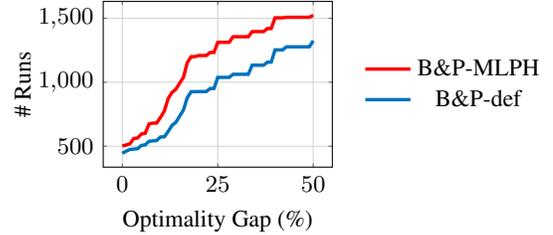
\begin{figure}[t!]
    \centering
	\begin{tikzpicture}
	\begin{axis}[height=0.45\columnwidth, width=0.55\columnwidth, grid style={line width=.1pt, draw=gray!10},major grid style={line width=.2pt,draw=gray!30}, xmajorgrids=true, ymajorgrids=true,  major tick length=0.05cm, minor tick length=0.0cm, ylabel = \small \# Runs, y label style={at={(axis description cs:-0.275,.5)},anchor=south}, xlabel = \small Optimality Gap (\%), xtick = {0, 25, 50}, ytick = {500, 1000, 1500}, legend style={at={(1.1,0.5)},anchor=west, font=\small,draw=none}, ticklabel style = {font=\small}]
	\addplot[color=red, line width=0.45mm] table [x=x, y=y, col sep=comma]
	{data/bp/gap-curve/mlhp-mix-10n-0.1n.txt};\addlegendentry{\small B\&P-MLPH}
	\addplot[color=blue, line width=0.45mm] table [x=x, y=y, col sep=comma] {data/bp/gap-curve/greedy.txt};\addlegendentry{\small B\&P-def}
    \end{axis} 
	\end{tikzpicture}
    \caption{The number of runs where a GCP instance can be solved within a certain optimality gap threshold.}
    \label{fig:bp_gap}
\end{figure}

Figure~\ref{fig:bp_gap} shows the number of runs for which B\&P-def and B\&P-MLPH can obtain optimality gap within a certain threshold value. We can observe that 1) B\&P-MLPH (red) solves GCP within a certain optimality gap in more runs than B\&P-def (blue); 2) B\&P-MLPH solves GCP to optimality ($Gap = 0\%$) in $482$ runs, better than $444$ runs by B\&P-def; 3) B\&P-MLPH obtains optimality gap (i.e., the LP at the root node is solved) on GCPs in $1524$ runs, whereas B\&P-def obtains optimality gap in $1325$ runs.

\begin{table}[t!]
    \centering
    \resizebox{0.85\linewidth}{!}{\begin{tabular}{@{}lrlrr@{}}
        \toprule
        \multirow{2}{*}{Instance} & \multirow{2}{*}{\# Nodes} &
        \multirow{2}{*}{Density} & 
        \multicolumn{2}{c}{Solving Time in Seconds} \\
         & & & B\&P-MLPH & B\&P-def\\
        \cmidrule(lr){1-3}\cmidrule(lr){4-5}
        
         r125.5 & 109 & 0.565 & \textbf{12} & 13 \\
         le450\_25b & 294 & 0.29 & \textbf{14} & 16 \\
         school1 & 355 & 0.603 & \textbf{51} & 1966 \\
         ash331GPIA & 661 & 0.038 & \textbf{53} & 307 \\
         qg.order30 & 900 & 0.129 & \textbf{58} & 60 \\
         will199GPIA & 660 & 0.054 & \textbf{70} & 128 \\
         flat300\_20\_0 & 300 & 0.953 & \textbf{84} & 343 \\
         DSJR500.1c & 311 & 0.972 & \textbf{99} & 109 \\
         flat300\_26\_0 & 300 & 0.965 & \textbf{223} & 1910 \\
        \cmidrule(lr){1-3}\cmidrule(lr){4-5}
        le450\_25a & 264 & 0.336 & 14 & \textbf{13} \\
         DSJC125.9 & 125 & 0.898 & 32 & \textbf{29} \\
         qg.order100 & 10000 & 0.04 & 6259 & \textbf{6127} \\
        \bottomrule
    \end{tabular}}
 \caption{Results for graphs solved by both B\&P methods in all runs.}
  \label{tab:bp1}
\end{table}

\begin{table}[t!]
    \centering
    \resizebox{0.85\linewidth}{!}{\begin{tabular}{@{}lrlrr@{}}
        \toprule
        \multirow{2}{*}{Instance} & \multirow{2}{*}{\# Nodes} &
        \multirow{2}{*}{Density} &
        \multicolumn{2}{c}{ Gap (\# Root Solved)} \\
         & & & B\&P-MLPH & B\&P-def \\
        \cmidrule(lr){1-3}\cmidrule(lr){4-5}
        
        queen16\_16 & 256 & 0.387 & \textbf{11.1} (\textbf{24}) & N/A (0) \\
        queen15\_15 & 225 & 0.411 & \textbf{11.8} (\textbf{24}) & N/A (0) \\
        le450\_25d & 433 & 0.366 & \textbf{10.7} (\textbf{24}) & N/A (0) \\
        le450\_15b & 410 & 0.187 & \textbf{6.2} (\textbf{24}) & N/A (0) \\
        le450\_25c & 435 & 0.362 & \textbf{10.7} (\textbf{24}) & N/A (0) \\
        le450\_15a & 407 & 0.189 & \textbf{6.2} (\textbf{24}) & N/A (0) \\
        DSJC250.9 & 250 & 0.896 & \textbf{2.3} (24)  & 4.1 (24)\\
        wap06a & 703 & 0.288 & \textbf{4.8} (\textbf{24}) & N/A (0) \\
        DSJC1000.9 & 1000 & 0.9 & 11.5 (\textbf{24}) & 11.5 (17) \\
        myciel6 & 95 & 0.338 & \textbf{31.0} (24)  & 42.9 (24)  \\
        qg.order40 & 1600 & 0.098 & \textbf{2.4} (\textbf{24}) & N/A (0) \\
        myciel5 & 47 & 0.437 & \textbf{22.9} (24) & 32.6 (24) \\
        
        \cmidrule(lr){1-3}\cmidrule(lr){4-5}
         1-Insertions\_5 & 202 & 0.121 & 43.8 (24) & \textbf{33.3} (24) \\
         2-Insertions\_5 & 597 & 0.044 & 50.0 (1) & 50.0 (\textbf{24})  \\
         3-Insertions\_5 & 1406 & 0.02 & N/A (0) & \textbf{50.0} (\textbf{24}) \\
         4-Insertions\_4 & 475 & 0.032 & 40.0 (15) & 40.0 (\textbf{24})  \\
         r1000.5 & 966 & 0.989 & 7.8 (24) & \textbf{1.2} (24) \\
        \bottomrule
    \end{tabular}}
     \caption{Results for graphs not solved by either of the two methods.}
  \label{tab:bp2}
\end{table}

\begin{table}[t!]
    \centering
    \resizebox{0.85\linewidth}{!}{\begin{tabular}{@{}lrlrr@{}}
        \toprule
        \multirow{2}{*}{Instance} & \multirow{2}{*}{\# Nodes} &
        \multirow{2}{*}{Density} &
        \multicolumn{2}{c}{\# Optimally Solved Runs} \\
         & & & B\&P-MLPH & B\&P-def \\
        \cmidrule(lr){1-3}\cmidrule(lr){4-5}
         le450\_5d & 450 & 0.193 & \textbf{21} & 10 \\
         le450\_5c & 450 & 0.194 & \textbf{21} & 1 \\
         ash608GPIA & 1215 & 0.021 & \textbf{24} & 0 \\
         2-Insertions\_3 & 37 & 0.216 & \textbf{17} & 9 \\
         DSJR500.5 & 486 & 0.972 & \textbf{17} & 15 \\
        \cmidrule(lr){1-3}\cmidrule(lr){4-5}
         1-FullIns\_4 & 38 & 0.364 & 18 & \textbf{23} \\
        queen9\_9 & 81 & 0.652 & 0 & \textbf{2} \\
        \bottomrule
    \end{tabular}}
 \caption{Results for graphs solved by the two methods in some~runs.}
  \label{tab:bp3}
\end{table}

Next, we report the numerical results on $36$ benchmark graphs where the performance of the two compared methods are significantly different (according to the student's $t$-test with a significance level of $0.05$). The results are grouped into Tables 4-6 based on their comparative performances. Table~\ref{tab:bp1} shows the results for graphs that can be optimally solved by both methods in all runs. Here, B\&P-MLPH uses less solving time than B\&P-def on $9$ graphs, and the speed-up is substantial on graphs such as `school1' ($38\times$) and `flat300\_26\_0' ($8\times$). In contrast, B\&P-def performs slightly better than B\&P-MLPH only on $3$ graphs. Table~\ref{tab:bp2} shows the results for hard graphs not solved by any method. B\&P-MLPH can still solve the LP at the root node for most graphs and runs. However, B\&P-def fails to solve the LP at the root node for many graphs, resulting in no optimality gap for those graphs. Table~\ref{tab:bp3} shows the results for the remaining graphs. B\&P-MLPH can solve GCP to optimality on more graphs and in more runs, compared to B\&P-def. Overall, B\&P-MLPH significantly outperforms B\&P-def on $26$ out of 36 graphs.
%  whereas B\&P-def performs slightly better than B\&P-MLPH on $10$ graphs
% Firstly, B\&P-MLPH can solve GCP to optimality on more graphs than B\&P-def. Secondly, B\&P-MLPH substantially outperforms B\&P-def as 1) B\&P-MLPH can solve GCP to optimality on more graphs than B\&P-def. 2) it solves GCP in much more runs than B\&P-def. 
% Overall, B\&P-MLPH outperforms B\&P-def on $26$ graphs, whereas B\&P-def performs slightly better on $10$ graphs.

Apart from these promising results, our studies also show that MLPH can be better integrated into B\&P based on certain conditions. Firstly, if the MWISP can be solved by an exact method efficiently on a graph, then it is not necessary to use MLPH (or any other pricing heuristic). In particular, on the set of $10$ graphs where B\&P-MLPH does not outperform B\&P-def, the performance of B\&P without using any pricing heuristic is on par with B\&P-def's. Secondly, if the improvement of the RMP becomes very slow, i.e., the tailing-off effect~\citep{gilmore1961linear}, an exact method can be used occasionally to solve MWISP to optimality even if MLPH can still find NRC columns. This is because MLPH is likely to keep finding NRC columns, which prevents the execution of the exact method from finding the optimal column and computing the Lagrangian lower bound. This reduces the chance of an early termination of CG. When applying this condition, we observe improved results of B\&P-MLPH on $12$ graphs. The detailed results can be found in the Appendix.

\section{Related Work}
Machine Learning for combinatorial optimization has received a lot of attention in recent years \citep{bengio2020machine}. Existing studies have applied ML in a variety of ways, such as learning variable selection methods \citep{khalil2016learning, gasse2019exact, liu2020learning, furian2021machine} or node selection methods \citep{he2014learning, furian2021machine} for exact branch-and-bound solvers; learning to select the best algorithm among its alternatives based on the problem characteristics~\citep{di2016dash,khalil2017learning}; learning to determine whether to perform problem reformulation~\citep{kruber2017learning, bonami2018learning} or problem reduction~\citep{sun2019using, ding2020accelerating}; learning primal heuristics aiming to construct an optimal solution directly \citep{dai2017learning,kool2018attention}; and learning to select columns for column generation~\citep{Morabit2020mlcs}.

% To our knowledge, we are the first to solve GCP using B\&P with a ML-based pricing heuristic. Our work is in line with the concurrent studies~\citep{quesnelWDS22,abs-2201-02535}, which tackle \emph{routing problems} using B\&P where ML is used to enhance a labeling algorithm for pricing. Our work is also related to recent studies in predicting optimal solutions for combinatorial optimization problems~\citep{li2018combinatorial,sun2020generalization, ding2020accelerating}. These solution-prediction-based methods focus on \emph{effectively finding a single best solution} for a single problem instance, whereas our pricing heuristic aims to \emph{efficiently generate many high-quality solutions} by solving a series of pricing problems.

Our work is in line with the recent studies in predicting optimal solutions for combinatorial optimization problems~\citep{li2018combinatorial,sun2019using,ding2020accelerating,sun2020generalization,sun2020boosting}. However, we are the first to leverage machine learning to develop a pricing heuristic for CG, to our knowledge. More specifically, the existing machine-learning-based heuristic methods focus on \emph{effectively finding a single best solution} (hopefully an optimal one) for a single problem instance, while our pricing heuristic aims to \emph{efficiently generate many high-quality solutions} by solving a series of pricing problems.

\section{Conclusion \& Future Work}

This paper presents a Machine-Learning-based Pricing Heuristic (MLPH) for tackling NP-hard pricing problems repeatedly encountered in the process of Column Generation~(CG). Specifically, we employ Support Vector Machine with linear kernel to fast predict `the optimal solution' for an NP-hard pricing problem, which is then adopted by a sampling-based method to \emph{construct many high-quality columns efficiently}. 

On the graph coloring problem, we demonstrate the efficacy of MLPH for solving its pricing problem - the maximum weight independent set problem. We demonstrate that MLPH can generate many more high-quality columns efficiently than existing state-of-the-art exact and heuristic methods. As a result, MLPH can significantly reduce the CG's computational time and enhance the branch-and-price exact method.

\added[id=YS]{In future work, we would like to extend our MLPH method to other combinatorial optimization problems such as vehicle routing problems. Our overarching aim is to develop a \emph{generic} ML-based pricing heuristic to speed up CG and branch-and-price for solving the Dantzig-Wolfe reformulation of combinatorial optimization problems.}

% \added[id=yunzhuang]{In future work, we would like to study whether MLPH can be adapted to handle different types of pricing problems efficiently. The key challenge will be obtaining effective ML prediction for sampling high-quality columns. We consider to make use of problem-specific knowledge in designing features, and we consider to explore deep neural models that learn feature representations automatically.}

% Since MLPH alleviates human efforts by ML, adapting MLPH to another problem can only require the use of problem-specific features and require the sampling method to properly handle constraints.
% Further, we would also examine the effector of different ML models in MLPH for solution~prediction.

% \section{Discussions}

% This paper presents a Machine-Learning-based Pricing Heuristic (MLPH) for tackling NP-hard pricing problems repeatedly encountered in the process of Column Generation~(CG). Specifically, we employ Support Vector Machine with linear kernel to fast predict `the optimal solution' for an NP-hard pricing problem, which is then adopted by a sampling-based method to \emph{construct many high-quality columns efficiently}. 

% On the graph coloring problem, we demonstrate the efficacy of MLPH for solving its pricing problem - the maximum weight independent set problem. We demonstrate that MLPH can generate many more high-quality columns efficiently than existing state-of-the-art exact and heuristic methods. As a result, MLPH can significantly reduce the CG's computational time and enhance the branch-and-price exact method.
% \bibliographystyle{aaai22}
\bibliography{lib}
\end{document}

\maketitle

This document is organized as follows. We present the column-selection algorithm, `replace-existing'; additional details about training our MLPH; a description about the B\&P implementation provided by SCIP~\citep{gamrath2020scip}; a discussion on how to adapt MLPH to other combinatorial optimization problems. Then, we present additional empirical results and analysis for CG and B\&P. For some abbreviations (used here and also in our main paper), their full terms can be found in Table~\ref{tab:abbv}.

\section{The Column Selection Method: `replace-existing'} 
The `replace-existing' method can keep a fixed size of the RMP throughout the process of CG, by replacing some of the columns already in the current RMP with some newly generated columns at each iteration of CG. We describe the algorithmic procedure as follows. We initialize an empty array to store the MISs (i.e., columns) selected for the next RMP. Each entry in this array is associated with a vertex in the graph, and an entry can store up to~$k$ MISs containing that vertex, where $k$ is the ratio between the number of MISs and the number of vertices in the graph. To ensure the solution of the next RMP is no worse than the current one, we firstly add into the array the MISs with non-zero values in the optimal solution of the current RMP. Then, we sort the rest of the MISs in ascending order of their reduced costs, and attempt to add them into the array one by one. A MIS can be added into an entry if satisfying two conditions: 1) the MIS contains the vertex associated with that entry and 2) the entry has less than~$k$ MISs. We note that the second condition is for maintaining the diversity of the columns in the next RMP, which is of importance~\cite{lubbecke2010column}. Lastly, we use the MISs in the array to form the RMP to start the next iteration.

\begin{table}[t]
    \centering
    \resizebox{.9\columnwidth}{!}{\begin{tabular}{ll}
    \toprule 
        Full Term & Abbreviation\\ \midrule
        
        Column Generation & CG \\
        Branch-and-Price & B\&P\\
        Machine-Learning-based Pricing Heuristic & MLPH \\
        Support Vector Machine with linear kernal & linear-SVM \\
        Graph Coloring Problem & GCP\\
        Linear-Programming relaxation & LP\\
        Restricted Master Problem & RMP \\
        Maximum Weight Independent Set Problem & MWISP\\
        Maximal Independent Set & MIS \\\bottomrule

    \end{tabular}}
    \caption{Terms and their abbreviations used in this document.}
    \label{tab:abbv}
\end{table}

\section{More Details of Training MLPH}

\begin{table*}[th!]
    \centering
    \resizebox{0.8\textwidth}{!}{\begin{tabular}{cccccc}
    \toprule 
        Feature & \begin{tabular}{@{}c@{}}Ranking-based \\ measure\end{tabular} & \begin{tabular}{@{}c@{}}Correlation-based \\ measure\end{tabular} & Vertex weight & Vertex degree & Vertex upper bound \\ \midrule
        
        Coefficient & 1.6557 & -1.0619 & -4.6320 & -1.5342 & 5.4064 \\\bottomrule
    \end{tabular}}
    \caption{The feature coefficients learned by linear-SVM. The value of the intercept term in linear-SVM is $1.1727$.}
    \label{tab:svm}
\end{table*}

In this part, we describe for MLPH the training graphs, the training time, and the learned parameters. We use ten small graphs randomly selected from the graph coloring benchmarks to train our MLPH method. These graphs are ``3-FullIns\_4", ``queen12\_12", ``1-Insertions\_6", ``mug88\_25", ``DSJC125.5", ``flat300\_20\_0", ``flat300\_26\_0", ``DSJC1000.9", ``DSJC250.1", and ``queen11\_11". We run CG on each of the graph and collect the training data during the CG process. More specifically, we use an exact method TSM~\citep{jiang2018two} to optimally solve the pricing problems, and collect solved MWISPs with optimal MISs every five iterations of CG to ensure the diversity of the training dataset. This results in a training dataset consisting of approximately $20,000$ training examples.

Using the training data, linear-SVM can be trained within one second, and its parameters are reported in Table~\ref{tab:svm}. Then, we optimize the Logistic model $p(d_i; \beta_0, \beta_1) = \frac{1}{1+e^{\beta_0 d_i + \beta_1}}$ using Bayesian Optimization (BO)~\citep{nogueira2014bayesian}. Using the ten graphs (listed above), BO runs CG-MLPH $300$ times with different sets of the parameters $\beta_0$ and $\beta_1$ in the Logistic model for solving the MWISP at the initial CG iteration, and identifies the set of parameters, $\beta_0=9.7750$ and $\beta_1=12.5564$, that yield best performance. The BO process can take several hours to complete. 

\section{The B\&P Implementation}

We use the B\&P code from SCIP~\citep{gamrath2020scip}. In the process of CG, the B\&P uses an efficient greedy search~\citep{mehrotra1996column} as the pricing heuristic and an exact method called $t$-clique for solving the pricing problems. Besides that, this B\&P implementation incorporates some specialized techniques for solving GCPs from previous studies~\citep{mehrotra1996column,malaguti2011exact}, such as a branching method that operates on the original compact formulation of GCP, a local search method that finds the initial primal solution of GCP (i.e., an upper bound) to be used as initial LP columns for RMP, an early branching method that may terminate CG early and safely by comparing Lagrangian lower bound of the current RMP with its objective value. In addition, the B\&P can also benefit from the functionalities provided by the SCIP solver, such as primal heuristics for finding better feasible solutions and column-management algorithms for efficiently re-optimizing RMPs during the CG process. 

It should be noted that we have made the following changes to the default B\&P: 1) We used the formula for calculating the Lagrangian lower bound in~\citep{malaguti2011exact} instead of its default one, in order to ensure the exactness of the B\&P, and 2) we do not enforce limitations on the number of simplex iterations (SCIP uses simplex methods to solve LP) and the maximum number of columns generated during the CG process, following previous work.

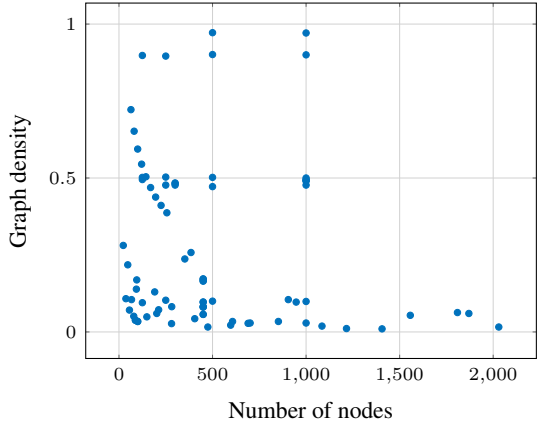
\begin{figure}[t!]
    \centering
    \begin{tikzpicture}
    	\begin{axis}[%
    	height=0.75\columnwidth, width=0.9\columnwidth, grid style={line width=.1pt, draw=gray!10},major grid style={line width=.2pt,draw=gray!30}, xlabel = \small Number of nodes, xtick = {0, 500, 1000, 1500, 2000}, xmajorgrids=true, ylabel= \small Graph density, ymajorgrids=true, major tick length=0.05cm, minor tick length=0.0cm, legend style={at={(1,1.2)}, column sep = 1pt, legend columns = 2,draw=none}, ticklabel style = {font=\scriptsize}%
    	]
    	\addplot[only marks,color=blue,%
    	    mark size=1.2pt]%
    	table[x=x, y=y,col sep=comma] {data/graph_stats_small.txt};
    	\end{axis}
    \end{tikzpicture}
    \caption{Characteristics of the set of $81$ small benchmark graphs.}
    \label{fig:gstats}
\end{figure}

\begin{table}[!htb]
\centering
   
        \resizebox{0.8\columnwidth}{!}{\begin{tabular}{lrlr}
             \toprule
            Graph & \# Nodes & Density & \begin{tabular}{@{}c@{}}Solving time of \\ the initial RMP\end{tabular}\\ \midrule

        % wap04a$^{*}$ & 5231 & 0.022 & 1800.0+ \\
        % wap03a$^{*}$ & 4730 & 0.026 & 1800.0+ \\
        % C4000.5$^{*}$ & 4000 & 0.500 & 1800.0+ \\
        % 4-FullIns\_5$^{*}$ & 4146 & 0.009 & 1280.5 \\
        % ash958GPIA$^{*}$ & 1916 & 0.007 & 938.8 \\
        % wap02a$^{*}$ & 2464 & 0.037 & 821.4 \\
        % C2000.5$^{*}$ & 2000 & 0.5 & 648.7 \\
        % wap01a$^{*}$ & 2368 & 0.04 & 582.6 \\
        wap08a & 1870 & 0.06 & 234.2 \\
        ash608GPIA & 1216 & 0.011 & 232.7 \\
        wap07a & 1809 & 0.063 & 230.9 \\
        abb313GPIA & 1557 & 0.054 & 161.7 \\
        3-FullIns\_5 & 2030 & 0.016 & 133.8 \\
        DSJC1000.1 & 1000 & 0.099 & 113.5 \\
        3-Insertions\_5 & 1406 & 0.01 & 98.8 \\
        r1000.1 & 1000 & 0.029 & 61.1 \\
        DSJC1000.5 & 1000 & 0.5 & 57.2 \\
        flat1000\_60\_0 & 1000 & 0.492 & 55.4 \\
        flat1000\_76\_0 & 1000 & 0.494 & 55.1 \\
        flat1000\_50\_0 & 1000 & 0.49 & 53.4 \\
        5-FullIns\_4 & 1085 & 0.019 & 30.2 \\
        will199GPIA & 701 & 0.029 & 22.4 \\
        wap05a & 905 & 0.105 & 21.8 \\
        wap06a & 947 & 0.097 & 21.3 \\
        DSJC1000.9 & 1000 & 0.9 & 10.2 \\
        DSJC500.1 & 500 & 0.1 & 8.4 \\
        2-FullIns\_5 & 852 & 0.034 & 7.7 \\
        4-Insertions\_4 & 475 & 0.016 & 6.9 \\
        2-Insertions\_5 & 597 & 0.022 & 6.1 \\
        4-FullIns\_4 & 690 & 0.028 & 5.6 \\
        r1000.5 & 1000 & 0.477 & 5.6 \\
        DSJC500.5 & 500 & 0.502 & 5.1 \\
        1-Insertions\_6 & 607 & 0.034 & 4.9 \\
        le450\_5a & 450 & 0.057 & 4.5 \\
        le450\_5b & 450 & 0.057 & 4.3 \\
        r1000.1c & 1000 & 0.971 & 3.0 \\
        le450\_25a & 450 & 0.082 & 2.8 \\
        le450\_15d & 450 & 0.166 & 2.6 \\
        le450\_15c & 450 & 0.165 & 2.6 \\
        le450\_15b & 450 & 0.081 & 2.3 \\
        le450\_15a & 450 & 0.081 & 2.3 \\
        le450\_25c & 450 & 0.172 & 2.1 \\
        le450\_25d & 450 & 0.172 & 2.1 \\
        le450\_5c & 450 & 0.097 & 1.8 \\
        queen16\_16 & 256 & 0.387 & 1.8 \\
        le450\_5d & 450 & 0.097 & 1.8 \\
        DSJC500.9 & 500 & 0.901 & 1.4 \\
        3-FullIns\_4 & 405 & 0.043 & 1.3 \\
        queen15\_15 & 225 & 0.411 & 1.2 \\
        3-Insertions\_4 & 281 & 0.027 & 1.2 \\
        school1 & 385 & 0.258 & 1.1 \\
        DSJR500.5 & 500 & 0.472 & 1.1 \\
        flat300\_26\_0 & 300 & 0.482 & 1.1 \\
        flat300\_28\_0 & 300 & 0.484 & 1.1 \\
        school1\_nsh & 352 & 0.237 & 1.1 \\
        flat300\_20\_0 & 300 & 0.477 & 1.0 \\
        DSJC250.1 & 250 & 0.103 & 0.8 \\
        queen14\_14 & 196 & 0.438 & 0.8 \\
        DSJC250.5 & 250 & 0.503 & 0.6 \\
        queen13\_13 & 169 & 0.469 & 0.5 \\
        DSJR500.1c & 500 & 0.972 & 0.5 \\
        1-FullIns\_5 & 282 & 0.082 & 0.4 \\
        1-Insertions\_5 & 202 & 0.06 & 0.4 \\
        queen12\_12 & 144 & 0.504 & 0.3 \\
        2-FullIns\_4 & 212 & 0.072 & 0.3 \\
        DSJC250.9 & 250 & 0.896 & 0.2 \\
        2-Insertions\_4 & 149 & 0.049 & 0.2 \\
        \bottomrule
        \end{tabular}}
        \caption{Statistics of the $81$ benchmark graphs labeled as `small'. The last attribute shows the solving time of the initial Restricted Master Problem (RMP) using the default LP method in Gurobi on a single CPU on our machine.}\label{tab:gstats1}
\end{table}

\begin{table}[!htb]
\centering
   
        \resizebox{0.8\columnwidth}{!}{\begin{tabular}{lrlr}
             \toprule
            Graph & \# Nodes & Density & \begin{tabular}{@{}c@{}}Solving time of \\ the initial RMP\end{tabular}\\ \midrule
        r250.5 & 250 & 0.477 & 0.2 \\
        queen11\_11 & 121 & 0.545 & 0.2 \\
        mug100\_1 & 100 & 0.034 & 0.2 \\
        mug100\_25 & 100 & 0.034 & 0.2 \\
        mug88\_1 & 88 & 0.038 & 0.1 \\
        myciel7 & 191 & 0.13 & 0.1 \\
        DSJC125.1 & 125 & 0.095 & 0.1 \\
        mug88\_25 & 88 & 0.038 & 0.1 \\
        queen10\_10 & 100 & 0.594 & 0.1 \\
        DSJC125.5 & 125 & 0.502 & 0.1 \\
        4-Insertions\_3 & 79 & 0.051 & 0.1 \\
        queen9\_9 & 81 & 0.652 & 0.1 \\
        r125.5 & 125 & 0.495 & 0.1 \\
        DSJC125.9 & 125 & 0.898 & 0.0 \\
        1-FullIns\_4 & 93 & 0.139 & 0.0 \\
        1-Insertions\_4 & 67 & 0.105 & 0.0 \\
        myciel6 & 95 & 0.169 & 0.0 \\
        3-Insertions\_3 & 56 & 0.071 & 0.0 \\
        queen8\_8 & 64 & 0.722 & 0.0 \\
        2-Insertions\_3 & 37 & 0.108 & 0.0 \\
        myciel5 & 47 & 0.218 & 0.0 \\
        myciel4 & 23 & 0.281 & 0.0 \\
        \bottomrule
        \end{tabular}}
        \caption{Statistics of the $81$ benchmark graphs labeled as `small' (Table~\ref{tab:gstats1} Continued).}\label{tab:gstats2}
\end{table}

\begin{table*}[ht!]
    \centering
    \resizebox{0.8\textwidth}{!}{\begin{tabular}{lccccccc}\toprule
     Methods & CG-MLPH &  CG-ACO & CG-Gurobi & CG-Gurobi-heur & CG-TSM & CG-Fastwclq & CG-LSCC \\ \midrule
      \# solved instances & \textbf{1096} & 916 & 895 & 728 & 932 & 964 & 812 \\
      LP objective value & \textbf{12.934} & 13.615 & 13.838 & 13.78 & 13.903 & 13.494 & 13.926 \\
      Solving time & \textbf{125.6} & 253.7 & 355.2 & 899.0 & 199.2 & 244.7 & 350.1 \\
        \bottomrule
    \end{tabular}}
    \caption{Solving statistics of CG using different pricing methods for $1944$ problem instances (generated using $81$ graphs labeled as `small'). Evaluation metrics include the total number of solved instances, the average objective value of the RMP (the lower the better), and the average solving time. The latter two statistics are obtained using the geometric mean over the results for all problem instances. }
    \label{tab:cg-solving}
\end{table*}

\begin{table*}[t!]
    \centering
    \resizebox{0.8\textwidth}{!}{\begin{tabular}{lrrrrrrr}\toprule
    Methods         & CG-MLPH & CG-ACO & CG-Gurobi & CG-Gurobi-heur & CG-TSM  & CG-Fastwclq & CG-LSCC \\ \midrule
      $1^{st}$ Place & \textbf{52} & 8 & 1 & 0 & 12 & 7 & 1 \\
      $2^{nd}$ Place \& Above & \textbf{69} & 23 & 10 & 13 & 26 & 20 & 1 \\
      $3^{rd}$ Place \& Above & \textbf{75} & 34 & 30 & 19 & 39 & 41 & 5 \\
    %   $4^{th}$ Place \& Above & 78 & 46 & 49 & 54 & 27 & 42 & 28 \\
        \bottomrule
    \end{tabular}}
    \caption{The number of graphs on which CG with a pricing method achieves top-$k$ ($k\in\{1,2,3\}$) performance among $81$ small graphs.}
    \label{tab:cg-rank}
\end{table*}

\begin{table*}[th!]
    \centering
    \resizebox{.75\textwidth}{!}{\begin{tabular}{@{}lrrrrrrr@{}}
        \toprule
        Graph & CG-MLPH & CG-ACO & CG-Gurobi & CG-Gurobi-heur& CG-TSM & CG-Fastwclq & CG-LSCC \\
        \cmidrule(lr){1-1}\cmidrule(lr){2-8}

        r1000.1 & \Large \textbf{2.1 (24)} & 2.7 (24) & 3.2 (24) & 2.8 (24) & 7.6 (11) & 6.0 (24) & 5.8 (24) \\
        will199GPIA & \Large \textbf{14.0 (23)} & N/A (0) & 78.8 (5) & 28.6 (5) & N/A (0) & N/A (0) & N/A (0) \\
        wap05a & \Large \textbf{2.1 (24)} & N/A (0) & 9.4 (24) & 6.7 (24) & 7.2 (24) & N/A (0) & N/A (0) \\
        wap06a & \Large \textbf{17.5 (12)} & N/A (0) & N/A (0) & N/A (0) & N/A (0) & N/A (0) & N/A (0) \\
        4-FullIns\_4 & \Large \textbf{18.2 (24)} & N/A (0) & 104.9 (24) & 37.4 (17) & N/A (0) & 26.6 (22) & N/A (0) \\
        le450\_25a & \Large \textbf{2.0 (24)} & 2.8 (24) & 3.7 (24) & 4.7 (24) & 3.0 (24) & 22.5 (24) & 38.2 (24) \\
        le450\_15b & \Large \textbf{10.0 (24)} & N/A (0) & 131.2 (24) & 22.7 (24) & N/A (0) & N/A (0) & N/A (0) \\
        le450\_15a & \Large \textbf{11.5 (24)} & N/A (0) & 139.8 (24) & 23.3 (24) & N/A (0) & N/A (0) & N/A (0) \\
        le450\_25c & \Large \textbf{20.6 (24)} & N/A (0) & N/A (0) & N/A (0) & N/A (0) & N/A (0) & N/A (0) \\
        le450\_25d & \Large \textbf{17.0 (24)} & N/A (0) & N/A (0) & N/A (0) & N/A (0) & N/A (0) & N/A (0) \\
        queen16\_16 & \Large \textbf{12.3 (24)} & 15.0 (24) & 86.4 (24) & 15.8 (24) & N/A (0) & 54.5 (17) & 418.5 (24) \\
        3-FullIns\_4 & \Large \textbf{8.4 (24)} & 16.8 (24) & 60.7 (24) & 23.9 (24) & 44.1 (24) & 12.2 (24) & 406.3 (23) \\
        queen15\_15 & \Large \textbf{10.4 (24)} & 12.6 (24) & 79.5 (24) & 12.8 (24) & N/A (0) & 44.8 (22) & 356.7 (24) \\
        queen14\_14 & \Large \textbf{10.0 (24)} & 10.9 (24) & 76.7 (24) & 11.1 (24) & 20.5 (24) & 51.1 (23) & 312.7 (24) \\
        queen13\_13 & \Large \textbf{7.5 (24)} & 7.9 (24) & 61.9 (24) & 8.0 (24) & 18.8 (24) & 31.9 (24) & 237.5 (24) \\
        1-FullIns\_5 & \Large \textbf{5.9 (24)} & 34.0 (1) & 534.4 (15) & N/A (0) & N/A (0) & 16.3 (24) & N/A (0) \\
        2-FullIns\_4 & \Large \textbf{4.8 (24)} & 7.0 (24) & 28.2 (24) & 12.5 (24) & 24.6 (24) & 5.7 (24) & 672.8 (24) \\
        DSJC250.9 & \Large \textbf{4.8 (24)} & 4.9 (24) & N/A (0) & N/A (0) & 18.7 (24) & 21.9 (24) & 66.0 (24) \\
        myciel6 & \Large \textbf{15.3 (24)} & 628.8 (24) & 92.9 (24) & 26.0 (24) & 85.0 (24) & 29.9 (24) & 934.4 (24) \\
        myciel5 & \Large \textbf{4.4 (24)} & 19.2 (24) & 14.3 (24) & 4.6 (24) & 12.7 (24) & 12.4 (24) & 74.6 (24) \\
        mug100\_1 & \Large \textbf{17.2 (24)} & 23.7 (24) & 50.8 (24) & N/A (0) & 29.0 (24) & 70.8 (24) & 214.3 (24) \\
        mug100\_25 & \Large \textbf{13.2 (24)} & 19.4 (24) & 45.0 (24) & 53.0 (2) & 24.8 (24) & 64.1 (24) & 193.5 (24) \\
        mug88\_1 & \Large \textbf{10.4 (24)} & 13.5 (24) & 33.2 (24) & 49.4 (14) & 18.4 (24) & 49.8 (24) & 120.2 (24) \\
        mug88\_25 & \Large \textbf{11.9 (24)} & 13.6 (24) & 38.1 (24) & 50.8 (5) & 21.6 (24) & 46.3 (24) & 148.3 (24) \\
        r125.5 & \Large \textbf{2.4 (24)} & \Large \textbf{2.4 (24)} & 6.1 (24) & 2.6 (24) & 4.4 (24) & 4.5 (24) & 12.8 (24) \\
        DSJC125.9 & \Large \textbf{2.9 (24)} & \Large \textbf{2.9 (24)} & 4.8 (24) & 3.8 (24) & 4.8 (24) & 5.3 (24) & 7.6 (24) \\
        1-FullIns\_4 & \Large \textbf{1.9 (24)} & \Large \textbf{1.9 (24)} & 14.8 (24) & 2.3 (24) & 15.0 (24) & 2.1 (24) & 23.2 (24) \\
        DSJR500.1c & \Large \textbf{2.0 (24)} & \Large \textbf{2.0 (24)} & 12.5 (11) & 20.8 (13) & 2.8 (24) & 2.5 (24) & 56.2 (24) \\
        myciel4 & \Large \textbf{2.0 (24)} & \Large \textbf{2.0 (24)} & 2.8 (24) & \Large \textbf{2.0 (24)} & 3.0 (24) & 2.8 (24) & 3.8 (24) \\
        r250.5 & 5.0 (24) & \Large \textbf{4.5 (24)} & 34.7 (24) & 7.1 (24) & 16.2 (24) & 115.9 (24) & 89.8 (24) \\
        DSJC125.5 & 283.8 (24) & \Large \textbf{16.2 (24)} & 51.6 (24) & 33.3 (6) & 40.3 (24) & 50.7 (24) & 412.4 (24) \\
        r1000.5 & 48.4 (14) & \Large \textbf{23.7 (24)} & N/A (0) & N/A (0) & 79.8 (24) & N/A (0) & N/A (0) \\
        r1000.1c & 4.0 (24) & \Large \textbf{3.9 (24)} & N/A (0) & N/A (0) & 20.7 (24) & 11.5 (24) & 623.0 (1) \\
        le450\_5c & N/A (0) & \Large \textbf{18.0 (7)} & N/A (0) & N/A (0) & N/A (0) & N/A (0) & N/A (0) \\
        le450\_5d & N/A (0) & \Large \textbf{21.8 (11)} & N/A (0) & N/A (0) & N/A (0) & N/A (0) & N/A (0) \\
        DSJC500.9 & 6.6 (24) & \Large \textbf{5.7 (24)} & N/A (0) & N/A (0) & 47.5 (24) & 63.8 (24) & 284.5 (24) \\
        DSJR500.5 & 7.3 (24) & \Large \textbf{5.0 (24)} & N/A (0) & N/A (0) & 28.1 (24) & 100.5 (21) & 196.0 (24) \\
        DSJC250.5 & 1025.0 (5) & \Large \textbf{65.7 (24)} & N/A (0) & N/A (0) & 96.3 (24) & 154.3 (24) & 1294.6 (17) \\
        DSJC1000.9 & 20.6 (24) & \Large \textbf{8.7 (24)} & N/A (0) & N/A (0) & 148.9 (24) & N/A (0) & N/A (0) \\
        queen12\_12 & 7.3 (24) & \Large \textbf{7.2 (24)} & 59.4 (24) & \Large \textbf{7.2 (24)} & 19.5 (24) & 32.0 (24) & 211.8 (24) \\
        queen8\_8 & 2.1 (24) & \Large \textbf{2.0 (24)} & 5.4 (24) & \Large \textbf{2.0 (24)} & 4.0 (24) & 4.0 (24) & 9.9 (24) \\    
        queen11\_11 & 7.0 (24) & 6.4 (24) & 53.6 (24) & \Large \textbf{6.0 (24)} & 21.3 (24) & 27.1 (24) & 176.4 (24) \\
        5-FullIns\_4 & 17.7 (10) & N/A (0) & N/A (0) & \Large \textbf{26.0 (18)} & N/A (0) & N/A (0) & N/A (0) \\
        4-Insertions\_3 & 1154.5 (24) & 422.8 (24) & 83.5 (24) & \Large \textbf{29.0 (24)} & 48.4 (24) & 58.2 (24) & 513.0 (24) \\
        queen9\_9 & 5.6 (24) & 4.5 (24) & 35.6 (24) & \Large \textbf{4.4 (24)} & 18.5 (24) & 21.1 (24) & 105.5 (24) \\
        1-Insertions\_4 & 31.7 (24) & 115.6 (24) & 42.5 (24) & \Large \textbf{16.0 (24)} & 34.6 (24) & 22.8 (24) & 311.7 (24) \\
        3-Insertions\_3 & 75.6 (24) & 73.0 (24) & 33.5 (24) & \Large \textbf{10.0 (24)} & 22.9 (24) & 28.1 (24) & 255.6 (24) \\
        2-Insertions\_3 & 7.6 (24) & 6.6 (24) & 10.6 (24) & \Large \textbf{3.8 (24)} & 9.4 (24) & 10.2 (24) & 27.9 (24) \\
        queen10\_10 & 6.5 (24) & 5.4 (24) & 52.8 (24) & \Large \textbf{5.2 (24)} & 25.2 (24) & 27.2 (24) & 175.8 (24) \\
        flat300\_28\_0 & N/A (0) & 111.2 (9) & N/A (0) & N/A (0) & \Large \textbf{124.5} (24) & 280.2 (19) & 1144.0 (3) \\
        DSJC125.1 & N/A (0) & N/A (0) & 172.7 (24) & N/A (0) & \Large \textbf{92.9 (24)} & N/A (0) & N/A (0) \\
        1-Insertions\_5 & N/A (0) & N/A (0) & N/A (0) & N/A (0) & N/A (0) & \Large \textbf{102.0 (24)} & N/A (0) \\
        2-Insertions\_4 & N/A (0) & N/A (0) & 308.7 (24) & N/A (0) & 225.2 (24) & \Large \textbf{105.2 (24)} & N/A (0) \\
        myciel7 & 178.3 (24) & N/A (0) & 412.5 (24) & N/A (0) & 465.6 (9) & \Large \textbf{57.2 (24)} & N/A (0) \\
        \bottomrule
    \end{tabular}}
    \caption{The number of CG iterations spent averaged over optimally solved test instances (out of $24$) and the number of optimally solved instances (in brackets). If a method cannot solve any instance for certain graphs, then `N/A' will be placed. Note that there are $54$ graphs whose problem instances can be optimally solved at least once by CG. It can be seen that MLPH helps CG capture optimal LP columns with a fewer number of iterations on $30$ of the graphs, and the performance of CG-MLPH is also competitive for the remaining graphs with a few exceptions.}
    \label{tab:cg-iter_pricing}
\end{table*}
\begin{table*}[t!]
    \centering
    \resizebox{0.75\textwidth}{!}{\begin{tabular}{lccccccc}\toprule
      Method       & MLPH & ACO & Gurobi & Gurobi-heur & TSM  & Fastwclq & LSCC \\ \midrule
      Average number of NRC columns found & \textbf{4017.7} & 807.5 & 4.5 & 474.9 & 78.0 & 90.9 & 0.9 \\
      The reduced cost of the best-found column & \textbf{-0.62} & -0.36 & -0.57 & -0.57 & -0.56 & -0.48 & 0.3 \\
        \bottomrule
    \end{tabular}}
    \caption{The performance of a pricing method for solving the MWISP at the initial iteration of the CG. The statistics are an average of the results of $1944$ problem instances using geometric mean.}
    \label{tab:cg-iter_pricing-stats}
\end{table*}

\section{Adapting MLPH to Other Problems}

In order to adapt MLPH to another combinatorial optimization problem, designing problem-specific features and modifying the sampling method will be required. Take the Vehicle Routing Problems as an example, where the pricing problem is the Resource-Constrained Shortest Path Problem (RCSPP). We can make use of the problem-specific features (e.g, the length and different types of resource consumption of an edge) and our statistical features, and train a linear-SVM to predict whether an edge is a part of the optimal path. In the testing phase, the prediction of linear-SVM can be used to guide the sampling of high-quality paths (i.e., columns).

% In future work, we would like to study whether MLPH can be adapted to handle different types of pricing problems efficiently. The key challenge will be obtaining effective ML prediction for sampling high-quality columns. We consider to make use of problem-specific knowledge in designing features, and we consider to explore deep neural models that learn feature representations automatically.

\section{Extended Results and Analysis for CG}  

Figure~\ref{fig:gstats} shows the set of $81$ benchmark graphs labeled as `small'. It can be seen that most graphs have less than $1000$ vertices and their densities are around or less than $0.5$. The rest of the graphs either are relatively dense or have a relatively large number of vertices. The statistics of individual graphs are shown in Table~\ref{tab:gstats1} and~\ref{tab:gstats2} for completeness. We generate $24$ instances on each of the benchmark graphs using different random seeds, resulting in a total number of $1944$ instances. Using these instances, we test CG with different heuristic-pricing methods under a cutoff time of $1800$ seconds, and the numerical results are reported as follows.

Table~\ref{tab:cg-solving} compares the performance of CG using different pricing methods.  Overall, our CG-MLPH consistently outperforms other methods by a large margin, in terms of all three metrics. With respect to the total number of solved instances, CG-MLPH can solve $1096$ instances as compared to $964$ instances solved by the second-best method, CG-Fastwclq. In terms of the average objective value of the RMP (indicating how close the current RMP solution is to the optimal LP solution), CG-MLPH achieves a value of $12.934$ whereas the second-best method CG-Fastwclq obtains a value of $13.494$. For solving time, CG-MLPH achieves an average of $125.6$ seconds, much better than $199.2$ seconds by the second-best method CG-TSM. Further, we report in Table~\ref{tab:cg-rank} that CG-MLPH achieves top-$k$ ($k\in\{1,2,3\}$) performances on individual graphs many more times than other methods. CG-ACO, CG-Gurobi, CG-TSM, and CG-Fastwclq are comparable with each other, but their performances are less competitive to CG-MLPH's.  

Table~\ref{tab:cg-iter_pricing} shows, for $54$ graphs whose problem instances can be solved by CG using different pricing methods at least once, the number of iterations spent (averaged over solved instances) by CG with each pricing method. As can be seen, CG-MLPH spends the fewest number of iterations for $30$ graphs (results highlighted in bold), and it is also very competitive for many graphs in the remaining $24$ graphs. In other words, MLPH helps CG better capture columns in an optimal LP solution using fewer CG iterations than other pricing methods typically. This is because MLPH can find \emph{many more high-quality columns} than other pricing methods, as evidenced by the results in Table~\ref{tab:cg-iter_pricing-stats}. Comparing to the sampling-based method ACO, the quality of the best-found solution by MLPH is much better. Comparing to more effective methods such as Gurobi or TSM, MLPH can find \emph{many more} NRC columns than those methods, and note that MLPH is still highly competitive in terms of the quality of the best-found column. Noticeably, CG-ACO and CG-Gourbi-heur are also competitive in terms of the iteration number (Table~\ref{tab:cg-iter_pricing}), because they can generate multiple NRC columns~(Table~\ref{tab:cg-iter_pricing-stats}).

Besides those summary statistics, we also include the results for individual graphs for reference: Table~\ref{tab:cg-solving1} and Table~\ref{tab:cg-solving2} show the results for CG using different pricing methods; Table~\ref{tab:cg-pp1} and Table~\ref{tab:cg-pp2} show the results for a pricing method for solving the MWISP at the initial CG iteration.

\section{Extended Results and Analysis for B\&P}

In this part, we firstly demonstrate that we can better integrate MLPH with B\&P based on some conditions. Then, we report parameter studies for B\&P-MLPH.

\begin{table*}[t!]
    \centering
    \resizebox{0.95\linewidth}{!}{\begin{tabular}{@{}lrlrrrrrrrr@{}}
        \toprule
         \multirow{2}{*}{Graph} & \multirow{2}{*}{\# Nodes} &
        \multirow{2}{*}{Density} &
        \multicolumn{2}{c}{\# optimally solved} &
        \multicolumn{2}{c}{\# root solved} & \multicolumn{2}{c}{\begin{tabular}{@{}c@{}}Optimality gap \\ (when the root node is solved) \end{tabular}} & 
        \multicolumn{2}{c}{\begin{tabular}{@{}c@{}}Total time\\ (in Seconds)\end{tabular}} \\
         & & & B\&P-def & B\&P-exact & B\&P-def & B\&P-exact & B\&P-def & B\&P-exact & B\&P-def & B\&P-exact \\
        \cmidrule(lr){1-3}\cmidrule(lr){4-5}\cmidrule(lr){6-7}\cmidrule(lr){8-9}\cmidrule(lr){10-11}
        1-FullIns\_4 & 38 & 0.364 & 23 & 24 & 24 & 24 & 0.8 & 0.0 & 980 & 81 \\
        queen9\_9 & 81 & 0.652 & 2 & 24 & 24 & 24 & 9.2 & 0.0 & 7517 & 19 \\
        DSJC125.9 & 125 & 0.898 & 24 & 24 & 24 & 24 & 0.0 & 0.0 & 29 & 23 \\
        r1000.5 & 966 & 0.989 & 0 & 0 & 24 & 24 & 1.2 & 8.7 & 8000 & 8000 \\
        \cmidrule(lr){1-3}\cmidrule(lr){4-5}\cmidrule(lr){6-7}\cmidrule(lr){8-9}\cmidrule(lr){10-11}
        le450\_25a & 264 & 0.336 & 24 & 24 & 24 & 24 & 0.0 & 0.0 & 13 & 15 \\
        qg.order100 & 10000 & 0.04 & 24 & 24 & 24 & 24 & 0.0 & 0.0 & 6127 & 6124 \\
        1-Insertions\_5 & 202 & 0.121 & 0 & 0 & 24 & 24 & 33.3 & 33.3 & 8000 & 8000 \\
        2-Insertions\_5 & 597 & 0.044 & 0 & 0 & 24 & 24 & 50.0 & 50.0 & 8000 & 8000 \\
        3-Insertions\_5 & 1406 & 0.02 & 0 & 0 & 24 & 24 & 50.0 & 50.0 & 8000 & 8000 \\
        4-Insertions\_4 & 475 & 0.032 & 0 & 0 & 24 & 24 & 40.0 & 40.0 & 8000 & 8000 \\
        \bottomrule
    \end{tabular}}
     \caption{Comparing B\&P using exact-pricing method directly (B\&P-exact) with B\&P-def that uses a greedy search on the set of $10$ graphs where B\&P-def is better than B\&P-MLPH. For the three graphs on the top, B\&P-exact outperforms B\&P-def substantially, whereas B\&P-def is better than B\&P-exact for the graph `r1000.5'. The two methods are comparable for the remaining graphs. The results indicate that the exact-pricing method is very efficient for solving MWISPs and hence, it is not necessary to use a heuristic-pricing method for these graphs. }
  \label{tab:bp-exact}
\end{table*}

\begin{table*}[t!]
    \centering
    \resizebox{0.95\linewidth}{!}{\begin{tabular}{@{}clrlrrrrrrrr@{}}
        \toprule
        \multirow{2}{*}{Group} & \multirow{2}{*}{Graph} & \multirow{2}{*}{\# Nodes} &
        \multirow{2}{*}{Density} &
        \multicolumn{2}{c}{\# optimally solved} &
        \multicolumn{2}{c}{\# root solved} & \multicolumn{2}{c}{\begin{tabular}{@{}c@{}}Optimality gap \\ (when the root node is solved) \end{tabular}} & 
        \multicolumn{2}{c}{\begin{tabular}{@{}c@{}}Total time\\ (in Seconds)\end{tabular}} \\
         & & & & B\&P-MLPH-def & B\&P-MLPH-forced & B\&P-MLPH-def & B\&P-MLPH-forced & B\&P-MLPH-def & B\&P-MLPH-forced & B\&P-MLPH-def & B\&P-MLPH-forced \\
        \cmidrule(lr){1-1}\cmidrule(lr){2-4}\cmidrule(lr){5-6}\cmidrule(lr){7-8}\cmidrule(lr){9-10}\cmidrule(lr){11-12}
        
         \multirow{14}{*}{Group 1} 
          & flat300\_26\_0 & 300 & 0.965 & 24 & 24 & 24 & 24 & 0.0 & 0.0 & 223 & 322 \\
         & qg.order30 & 900 & 0.129 & 24 & 24 & 24 & 24 & 0.0 & 0.0 & 58 & 61 \\
         & le450\_5c & 450 & 0.194 & 21 & 8 & 21 & 8 & 0.0 & 0.0 & 4521 & 6325 \\
         & qg.order40 & 1600 & 0.098 & 0 & 0 & 24 & 21 & 2.4 & 2.4 & 8002 & 8002 \\
         & will199GPIA & 660 & 0.054 & 24 & 24 & 24 & 24 & 0.0 & 0.0 & 70 & 76 \\
         & school1\_nsh & 326 & 0.547 & 24 & 24 & 24 & 24 & 0.0 & 0.0 & 626 & 1355 \\
         & le450\_25a & 264 & 0.336 & 24 & 24 & 24 & 24 & 0.0 & 0.0 & 13 & 14 \\
         & qg.order60 & 3600 & 0.066 & 24 & 24 & 24 & 24 & 0.0 & 0.0 & 631 & 658 \\
         & le450\_25b & 294 & 0.29 & 24 & 24 & 24 & 24 & 0.0 & 0.0 & 14 & 15 \\
         & wap05a & 665 & 0.314 & 24 & 24 & 24 & 24 & 0.0 & 0.0 & 74 & 81 \\
         & ash608GPIA & 1215 & 0.021 & 24 & 0 & 24 & 0 & 0.0 & N/A & 4899 & 8000 \\
         & r125.5 & 109 & 0.565 & 24 & 24 & 24 & 24 & 0.0 & 0.0 & 12 & 13 \\
         \cmidrule(lr){1-1}\cmidrule(lr){2-4}\cmidrule(lr){5-6}\cmidrule(lr){7-8}\cmidrule(lr){9-10}\cmidrule(lr){11-12}
        \multirow{14}{*}{Group 2} 
         & 2-Insertions\_3 & 37 & 0.216 & 17 & 20 & 24 & 24 & 7.3 & 4.2 & 5573 & 5472 \\
         & le450\_25d & 433 & 0.366 & 0 & 3 & 24 & 24 & 10.7 & 9.4 & 8000 & 7005 \\
         & queen9\_9 & 81 & 0.652 & 0 & 1 & 24 & 24 & 10.0 & 9.6 & 8000 & 7667 \\
         & queen16\_16 & 256 & 0.387 & 0 & 1 & 24 & 1 & 11.1 & 0.0 & 8000 & 7668 \\
         & 1-FullIns\_4 & 38 & 0.364 & 18 & 21 & 24 & 24 & 5.0 & 2.5 & 2967 & 1635 \\
         & 4-Insertions\_4 & 475 & 0.032 & 0 & 0 & 15 & 24 & 40.0 & 40.0 & 8001 & 8000 \\
         & DSJC125.5 & 125 & 0.502 & 0 & 0 & 24 & 24 & 15.4 & 13.7 & 8000 & 8000 \\
         & DSJC125.1 & 125 & 0.19 & 0 & 8 & 24 & 24 & 16.7 & 11.1 & 8000 & 5334 \\
         & queen15\_15 & 225 & 0.411 & 0 & 1 & 24 & 24 & 11.8 & 11.3 & 8000 & 7667 \\
         & 2-Insertions\_5 & 597 & 0.044 & 0 & 0 & 1 & 23 & 50.0 & 50.0 & 8000 & 8000 \\
         & le450\_5d & 450 & 0.193 & 21 & 24 & 23 & 24 & 1.4 & 0.0 & 3542 & 931 \\
         & le450\_25c & 435 & 0.362 & 0 & 1 & 24 & 24 & 10.7 & 10.3 & 8000 & 7668 \\
        \bottomrule
    \end{tabular}}
     \caption{Comparing the default setting of B\&P-MLPH (used in our main paper) with a variant that forces the exact method to run even if columns with negative reduced costs are found by MLPH, named B\&P-MLPH-forced. The exact method runs occasionally under a frequency of $5$ CG iterations when the objective value of the RMP improves less than $0.05$ in successive $5$ iterations. The results show that this variant is better than the default B\&P-MLPH on $14$ graphs (Group 2) and worse on a different set of $14$ graphs (Group 1). The results indicate that, for graphs in Group 2, finding the optimal solution by an exact method and computing the Lagrangian lower bound can benefit B\&P-MLPH from the early termination of CG. However, the exact method should be given a certain cutoff time, because it can take too much computational budget and degrade the performance of B\&P-MLPH for graphs in Group 1.}
  \label{tab:bp-tailoff}
\end{table*}

\begin{table*}[t!]
    \centering
    \resizebox{0.95\linewidth}{!}{\begin{tabular}{@{}clrlrrrrrrrr@{}}
        \toprule
        \multirow{2}{*}{Group} & \multirow{2}{*}{Graph} & \multirow{2}{*}{\# Nodes} &
        \multirow{2}{*}{Density} &
        \multicolumn{2}{c}{\# optimally solved} &
        \multicolumn{2}{c}{\# root solved} & \multicolumn{2}{c}{\begin{tabular}{@{}c@{}}Optimality gap \\ (when the root node is solved) \end{tabular}} & 
        \multicolumn{2}{c}{\begin{tabular}{@{}c@{}}Total time\\ (in Seconds)\end{tabular}} \\
         & & & & B\&P-MLPH-def & B\&P-MLPH-uni & B\&P-MLPH-def & B\&P-MLPH-uni & B\&P-MLPH-def & B\&P-MLPH-uni & B\&P-MLPH-def & B\&P-MLPH-uni \\
        \cmidrule(lr){1-1}\cmidrule(lr){2-4}\cmidrule(lr){5-6}\cmidrule(lr){7-8}\cmidrule(lr){9-10}\cmidrule(lr){11-12}
        
        \multirow{6}{*}{Group 1}
        & r125.5 & 109 & 0.565 & 24 & 24 & 24 & 24 & 0.0 & 0.0 & 12 & 13 \\
        & DSJR500.1c & 311 & 0.972 & 24 & 24 & 24 & 24 & 0.0 & 0.0 & 99 & 102 \\
        & 1-FullIns\_4 & 38 & 0.364 & 18 & 14 & 24 & 24 & 5.0 & 8.3 & 2967 & 4893 \\
        & le450\_5d & 450 & 0.193 & 21 & 19 & 23 & 24 & 1.4 & 3.5 & 3542 & 4229 \\
        & 2-Insertions\_3 & 37 & 0.216 & 17 & 7 & 24 & 24 & 7.3 & 17.7 & 5573 & 7317 \\
        & 4-Insertions\_4 & 475 & 0.032 & 0 & 0 & 15 & 14 & 40.0 & 40.0 & 8001 & 8017 \\

         \cmidrule(lr){1-1}\cmidrule(lr){2-4}\cmidrule(lr){5-6}\cmidrule(lr){7-8}\cmidrule(lr){9-10}\cmidrule(lr){11-12}
         \multirow{3}{*}{Group 2}
         & DSJC125.9 & 125 & 0.898 & 24 & 24 & 24 & 24 & 0.0 & 0.0 & 32 & 30 \\
         & DSJR500.5 & 486 & 0.972 & 17 & 21 & 24 & 24 & 0.2 & 0.1 & 3304 & 2055 \\
         & queen11\_11 & 121 & 0.545 & 0 & 1 & 24 & 24 & 15.4 & 14.7 & 8000 & 7737 \\
        \bottomrule
    \end{tabular}}
     \caption{Comparing the default setting of B\&P-MLPH (used in our main paper) with a variant that uses a uniform column-limit parameter $\theta=n$ for column selection, named B\&P-MLPH-uni. Note that the default B\&P-MLPH sets different column-limit parameters for the root node ($\theta=n$) and other nodes ($\theta=n$). The results show that the default B\&P-MLPH is significantly better than the compared variant on $6$ graphs (Group 1); the compared one is better on $3$ graphs (Group 2); the performances of the two methods are comparable for other graphs not presented here.}
  \label{tab:bp-columnlimit}
\end{table*}

Table~\ref{tab:bp-exact} shows that, on the set of $10$ graphs where B\&P-MLPH does not outperform B\&P-def, B\&P without using any heuristic-pricing method (i.e, only using the exact-pricing method) can perform on par to or even better than B\&P-def. In other words, the exact method is very efficient for solving MWISPs for these graphs. In such cases, it is not necessary to use a heuristic-pricing method (e.g., greedy search or MLPH). In practice, one can detect such graphs by running an exact method (with a small cutoff time) for solving the MWISP at the initial iteration of CG. Whether to use a pricing-heuristic method can then be determined based on the performance of the exact method on this graph. 

Table~\ref{tab:bp-tailoff} shows that, for some benchmark graphs, B\&P-MLPH can benefit from forcing the exact-pricing method to optimally solve MWISPs occasionally in the latter CG iterations. When enforcing this rule to B\&P-MLPH, the performance of B\&P-MLPH is significantly improved for a set of $14$ benchmark graphs. However, its performance is also significantly degraded for another set of $14$ graphs. These results indicate that, for certain graphs, finding the optimal solution by an exact method and computing the Lagrangian lower bound can early terminate CG. However, the exact method should be given a certain cutoff time, because it can introduce computational overhead and degrade the performance of B\&P-MLPH for some graphs.

\begin{table*}[t!]
    \centering
    \resizebox{0.95\linewidth}{!}{\begin{tabular}{@{}clrlrrrrrrrr@{}}
        \toprule
        \multirow{2}{*}{Group} & \multirow{2}{*}{Graph} & \multirow{2}{*}{\# Nodes} &
        \multirow{2}{*}{Density} &
        \multicolumn{2}{c}{\# optimally solved} &
        \multicolumn{2}{c}{\# root solved} & \multicolumn{2}{c}{\begin{tabular}{@{}c@{}}Optimality gap \\ (when the root node is solved) \end{tabular}} & 
        \multicolumn{2}{c}{\begin{tabular}{@{}c@{}}Total time\\ (in Seconds)\end{tabular}} \\
         & & & & B\&P-MLPH-$10n$ & B\&P-MLPH-$n$ & B\&P-MLPH-$10n$ & B\&P-MLPH-$n$ & B\&P-MLPH-$10n$ & B\&P-MLPH-$n$ & B\&P-MLPH-$10n$ & B\&P-MLPH-$n$ \\
        \cmidrule(lr){1-1}\cmidrule(lr){2-4}\cmidrule(lr){5-6}\cmidrule(lr){7-8}\cmidrule(lr){9-10}\cmidrule(lr){11-12}
        
        \multirow{13}{*}{Group 1} 
         & flat300\_26\_0 & 300 & 0.965 & 24 & 24 & 24 & 24 & 0.0 & 0.0 & 223 & 1200 \\
         & qg.order30 & 900 & 0.129 & 24 & 24 & 24 & 24 & 0.0 & 0.0 & 58 & 62 \\
         & le450\_5c & 450 & 0.194 & 21 & 4 & 21 & 4 & 0.0 & 0.0 & 4521 & 7301 \\
         & flat300\_20\_0 & 300 & 0.953 & 24 & 24 & 24 & 24 & 0.0 & 0.0 & 84 & 177 \\
         & r1000.5 & 966 & 0.989 & 0 & 0 & 24 & 24 & 7.8 & 9.3 & 8000 & 8000 \\
         & queen16\_16 & 256 & 0.387 & 0 & 0 & 24 & 0 & 11.1 & N/A & 8000 & 8000 \\
         & qg.order60 & 3600 & 0.066 & 24 & 24 & 24 & 24 & 0.0 & 0.0 & 631 & 658 \\
         & wap05a & 665 & 0.314 & 24 & 24 & 24 & 24 & 0.0 & 0.0 & 74 & 81 \\
         & ash608GPIA & 1215 & 0.021 & 24 & 18 & 24 & 18 & 0.0 & 0.0 & 4899 & 7403 \\
         & queen15\_15 & 225 & 0.411 & 0 & 0 & 24 & 19 & 11.8 & 11.8 & 8000 & 8000 \\
         & r125.5 & 109 & 0.565 & 24 & 24 & 24 & 24 & 0.0 & 0.0 & 12 & 14 \\
         & le450\_5d & 450 & 0.193 & 21 & 16 & 23 & 20 & 1.4 & 3.3 & 3542 & 4333 \\
         & qg.order100 & 10000 & 0.04 & 24 & 24 & 24 & 24 & 0.0 & 0.0 & 6259 & 6463 \\
        \cmidrule(lr){1-1}\cmidrule(lr){2-4}\cmidrule(lr){5-6}\cmidrule(lr){7-8}\cmidrule(lr){9-10}\cmidrule(lr){11-12}
         \multirow{15}{*}{Group 2} 
         & myciel6 & 95 & 0.338 & 0 & 0 & 24 & 24 & 31.0 & 28.6 & 8000 & 8000 \\
         & 1-Insertions\_4 & 67 & 0.21 & 0 & 0 & 24 & 24 & 37.5 & 20.8 & 8003 & 8000 \\
         & DSJC125.9 & 125 & 0.898 & 24 & 24 & 24 & 24 & 0.0 & 0.0 & 32 & 26 \\
         & 2-Insertions\_3 & 37 & 0.216 & 17 & 22 & 24 & 24 & 7.3 & 2.1 & 5573 & 3137 \\
         & will199GPIA & 660 & 0.054 & 24 & 24 & 24 & 24 & 0.0 & 0.0 & 70 & 61 \\
         & DSJR500.5 & 486 & 0.972 & 17 & 23 & 24 & 24 & 0.2 & 0.0 & 3304 & 1175 \\
         & myciel5 & 47 & 0.437 & 0 & 0 & 24 & 24 & 22.9 & 16.7 & 8007 & 8002 \\
         & queen9\_9 & 81 & 0.652 & 0 & 24 & 24 & 24 & 10.0 & 0.0 & 8000 & 18 \\
         & ash331GPIA & 661 & 0.038 & 24 & 24 & 24 & 24 & 0.0 & 0.0 & 53 & 39 \\
         & 1-FullIns\_4 & 38 & 0.364 & 18 & 23 & 24 & 23 & 5.0 & 0.0 & 2967 & 597 \\
         & 4-Insertions\_4 & 475 & 0.032 & 0 & 0 & 15 & 23 & 40.0 & 40.0 & 8001 & 8000 \\
         & DSJC125.5 & 125 & 0.502 & 0 & 1 & 24 & 24 & 15.4 & 10.1 & 8000 & 7932 \\
         & 1-FullIns\_5 & 78 & 0.277 & 0 & 0 & 24 & 24 & 31.9 & 16.7 & 8001 & 8001 \\
         & 2-Insertions\_5 & 597 & 0.044 & 0 & 0 & 1 & 14 & 50.0 & 50.0 & 8000 & 8000 \\
         & 1-Insertions\_5 & 202 & 0.121 & 0 & 0 & 24 & 24 & 43.8 & 36.8 & 8000 & 8000 \\
        \bottomrule
    \end{tabular}}
     \caption{Comparing B\&P-MLPH using sample size $10n$ (used in the main paper) and the one using sample size $n$, named B\&P-MLPH-$n$. Group 1 shows the graphs where using sample size $10n$ is better and group 2 shows the graphs where using sample size $n$ is better. The performances of the two methods are comparable to the rest of the benchmark graphs. The results show that no single sample size can fit all graph benchmarks. }
  \label{tab:bp-samplesize1}
\end{table*}

\begin{table*}[t!]
    \centering
    \resizebox{0.95\linewidth}{!}{\begin{tabular}{@{}clrlrrrrrrrr@{}}
        \toprule
        \multirow{2}{*}{Group} & \multirow{2}{*}{Graph} & \multirow{2}{*}{\# Nodes} &
        \multirow{2}{*}{Density} &
        \multicolumn{2}{c}{\# optimally solved} &
        \multicolumn{2}{c}{\# root solved} & \multicolumn{2}{c}{\begin{tabular}{@{}c@{}}Optimality gap \\ (when the root node is solved) \end{tabular}} & 
        \multicolumn{2}{c}{\begin{tabular}{@{}c@{}}Total time\\ (in Seconds)\end{tabular}} \\
         & & & & B\&P-MLPH-$10n$ & B\&P-MLPH-$0.1n$ & B\&P-MLPH-$10n$ & B\&P-MLPH-$0.1n$ & B\&P-MLPH-$10n$ & B\&P-MLPH-$0.1n$ & B\&P-MLPH-$10n$ & B\&P-MLPH-$0.1n$ \\
        \cmidrule(lr){1-1}\cmidrule(lr){2-4}\cmidrule(lr){5-6}\cmidrule(lr){7-8}\cmidrule(lr){9-10}\cmidrule(lr){11-12}
        
        \multirow{15}{*}{Group 1} 
         & flat300\_26\_0 & 300 & 0.965 & 24 & 24 & 24 & 24 & 0.0 & 0.0 & 223 & 1919 \\
         & le450\_5c & 450 & 0.194 & 21 & 0 & 21 & 0 & 0.0 & N/A & 4521 & 8000 \\
         & flat300\_20\_0 & 300 & 0.953 & 24 & 24 & 24 & 24 & 0.0 & 0.0 & 84 & 363 \\
         & school1 & 355 & 0.603 & 24 & 24 & 24 & 24 & 0.0 & 0.0 & 51 & 507 \\
         & r1000.5 & 966 & 0.989 & 0 & 0 & 24 & 24 & 7.8 & 9.3 & 8000 & 8000 \\
         & le450\_25d & 433 & 0.366 & 0 & 0 & 24 & 0 & 10.7 & N/A & 8000 & 8000 \\
         & le450\_15a & 407 & 0.189 & 0 & 0 & 24 & 1 & 6.2 & 6.2 & 8000 & 8000 \\
         & school1\_nsh & 326 & 0.547 & 24 & 20 & 24 & 20 & 0.0 & 0.0 & 626 & 3292 \\
         & queen16\_16 & 256 & 0.387 & 0 & 0 & 24 & 0 & 11.1 & N/A & 8000 & 8000 \\
         & wap06a & 703 & 0.288 & 0 & 0 & 24 & 1 & 4.8 & 4.8 & 8000 & 8000 \\
         & ash608GPIA & 1215 & 0.021 & 24 & 0 & 24 & 0 & 0.0 & N/A & 4899 & 8000 \\
         & queen15\_15 & 225 & 0.411 & 0 & 0 & 24 & 0 & 11.8 & N/A & 8000 & 8000 \\
         & r125.5 & 109 & 0.565 & 24 & 24 & 24 & 24 & 0.0 & 0.0 & 12 & 12 \\
         & le450\_25c & 435 & 0.362 & 0 & 0 & 24 & 0 & 10.7 & N/A & 8000 & 8000 \\
         & le450\_15b & 410 & 0.187 & 0 & 0 & 24 & 0 & 6.2 & N/A & 8000 & 8000 \\
        \cmidrule(lr){1-1}\cmidrule(lr){2-4}\cmidrule(lr){5-6}\cmidrule(lr){7-8}\cmidrule(lr){9-10}\cmidrule(lr){11-12}
         \multirow{20}{*}{Group 2} 
         & myciel6 & 95 & 0.338 & 0 & 0 & 24 & 24 & 31.0 & 28.6 & 8000 & 8000 \\
         & qg.order30 & 900 & 0.129 & 24 & 24 & 24 & 24 & 0.0 & 0.0 & 58 & 56 \\
         & 1-Insertions\_4 & 67 & 0.21 & 0 & 0 & 24 & 24 & 37.5 & 20.8 & 8003 & 8000 \\
         & DSJC125.9 & 125 & 0.898 & 24 & 24 & 24 & 24 & 0.0 & 0.0 & 32 & 25 \\
         & 2-Insertions\_3 & 37 & 0.216 & 17 & 20 & 24 & 24 & 7.3 & 4.2 & 5573 & 4451 \\
         & will199GPIA & 660 & 0.054 & 24 & 24 & 24 & 24 & 0.0 & 0.0 & 70 & 49 \\
         & DSJR500.5 & 486 & 0.972 & 17 & 24 & 24 & 24 & 0.2 & 0.0 & 3304 & 1653 \\
         & myciel5 & 47 & 0.437 & 0 & 0 & 24 & 24 & 22.9 & 16.7 & 8007 & 8000 \\
         & queen9\_9 & 81 & 0.652 & 0 & 24 & 24 & 24 & 10.0 & 0.0 & 8000 & 20 \\
         & ash331GPIA & 661 & 0.038 & 24 & 24 & 24 & 24 & 0.0 & 0.0 & 53 & 38 \\
         & le450\_25a & 264 & 0.336 & 24 & 24 & 24 & 24 & 0.0 & 0.0 & 13 & 12 \\
         & qg.order60 & 3600 & 0.066 & 24 & 24 & 24 & 24 & 0.0 & 0.0 & 631 & 614 \\
         & 1-FullIns\_4 & 38 & 0.364 & 18 & 23 & 24 & 23 & 5.0 & 0.0 & 2967 & 421 \\
         & 4-Insertions\_4 & 475 & 0.032 & 0 & 0 & 15 & 23 & 40.0 & 40.0 & 8001 & 8000 \\
         & le450\_25b & 294 & 0.29 & 24 & 24 & 24 & 24 & 0.0 & 0.0 & 14 & 13 \\
         & wap05a & 665 & 0.314 & 24 & 24 & 24 & 24 & 0.0 & 0.0 & 74 & 71 \\
         & 2-Insertions\_5 & 597 & 0.044 & 0 & 0 & 1 & 22 & 50.0 & 50.0 & 8000 & 8000 \\
         & le450\_5d & 450 & 0.193 & 21 & 22 & 23 & 23 & 1.4 & 0.7 & 3542 & 2005 \\
         & 1-Insertions\_5 & 202 & 0.121 & 0 & 0 & 24 & 24 & 43.8 & 34.7 & 8000 & 8000 \\
         & qg.order100 & 10000 & 0.04 & 24 & 24 & 24 & 24 & 0.0 & 0.0 & 6259 & 6125 \\
        \bottomrule
    \end{tabular}}
     \caption{Comparing B\&P-MLPH using sample size $10n$ (used in the main paper) and the one using sample size $0.1n$, named B\&P-MLPH-$0.1n$. Group 1 shows the graphs where using sample size $10n$ is better and group 2 shows the graphs where using sample size $0.1n$ is better. The results show that no single sample size can fit all graph benchmarks. }
  \label{tab:bp-samplesize2}
\end{table*}

Next, we report parameter studies for B\&P-MLPH. Table~\ref{tab:bp-columnlimit} compares the default setting of B\&P-MLPH (used in our main paper) with a variant that uses a constant column-limit parameter $\theta=n$ for column selection. Note that our default B\&P-MLPH sets column-limit parameter $\theta=n$ for the root node and $\theta=0.1n$ for other nodes. It can be seen that the performance of B\&P-MLPH with the default parameter setting is better overall. The reason is of two folds. Firstly, the child nodes often have a sufficient number of quality columns inherited from their parents, and hence keep adding a large number of columns can be less useful. Secondly, setting a small column limit parameter for child nodes can slow down the growth in the size of RMPs and reduce the computational overhead for an LP solver to solve RMPs. 

Table~\ref{tab:bp-samplesize1} and Table~\ref{tab:bp-samplesize2} present the results for varying the sampling size $\{10n,n,0.1n\}$ in B\&P-MLPH, which shows that no single sample size of MLPH fits all benchmark graphs. These results indicate that B\&P-MLPH can be possibly further improved by configuring B\&P-MLPH according to the characteristics of the problem at hand. More specifically, the sample size can be set adaptively using information collected from the solving process, such as the speed of the CG convergence when using different samples sizes. Further, selecting a suitable sample size for B\&P-MLPH can be viewed as an algorithm selection problem~\cite{di2016dash, khalil2016learning}. We leave these possible solutions to future work as building an advanced parameter selection model requires significant effort and is out of the scope of the current study.

\begin{landscape}
\begin{table}[h!]
    \centering
    \resizebox{\linewidth}{!}{\begin{tabular}{@{}lrrrrrrr|ccccccc|rrrrrrr@{}}
        \toprule
        \multirow{2}{*}{Graph} & \multicolumn{7}{c}{\# solve runs} & \multicolumn{7}{c}{LP objective value} & \multicolumn{7}{c}{Time} \\
         & MLPH & ACO & Gurobi & Gurobi-heur& TSM & Fastwclq & LSCC & MLPH & ACO & Gurobi & Gurobi-heur & TSM & Fastwclq & LSCC & MLPH & ACO & Gurobi & Gurobi-heur & TSM & Fastwclq & LSCC\\ 
        \cmidrule(lr){1-1}\cmidrule(lr){2-8}\cmidrule(lr){9-15}\cmidrule(lr){16-22}
        ash608GPIA & 0 & 0 & 0 & 0 & 0 & 0 & 0 & 3.372 & 3.368 & 3.394 & 3.392 & 3.394 & 3.385 & 3.394 & 2045.9 & 1988.8 & 1927.9 & 1887.1 & 1935.3 & 1903.0 & 1963.3 \\
        wap08a & 0 & 0 & 0 & 0 & 0 & 0 & 0 & 47.53 & 62.258 & 58.627 & 57.94 & 65.334 & 65.456 & 64.801 & 2086.6 & 1926.8 & 1906.5 & 1911.9 & 1961.7 & 2043.5 & 1884.0 \\
        wap07a & 0 & 0 & 0 & 0 & 0 & 0 & 0 & 47.105 & 60.4 & 59.423 & 57.184 & 63.366 & 63.496 & 62.555 & 2217.2 & 1959.4 & 2004.6 & 1983.0 & 2057.1 & 1987.6 & 1894.3 \\
        abb313GPIA & 0 & 0 & 0 & 0 & 0 & 0 & 0 & 8.455 & 9.371 & 9.345 & 9.334 & 9.54 & 9.483 & 9.588 & 2164.1 & 1902.6 & 1896.7 & 1904.0 & 1900.3 & 1908.6 & 2154.6 \\
        3-FullIns\_5 & 0 & 0 & 0 & 0 & 0 & 0 & 0 & 5.697 & 8.763 & 7.684 & 6.648 & 9.672 & 6.607 & 9.349 & 1903.7 & 1866.0 & 1865.7 & 1870.7 & 1866.6 & 1913.2 & 1840.9 \\
        DSJC1000.1 & 0 & 0 & 0 & 0 & 0 & 0 & 0 & 19.984 & 21.028 & 22.107 & 22.124 & 22.296 & 20.837 & 21.645 & 1890.5 & 1880.8 & 1863.8 & 1866.9 & 1977.2 & 1899.9 & 1851.4 \\
        3-Insertions\_5 & 0 & 0 & 0 & 0 & 0 & 0 & 0 & 3.343 & 4.425 & 3.974 & 3.622 & 4.609 & 3.866 & 4.871 & 1939.7 & 1858.7 & 1852.9 & 1899.9 & 1897.8 & 1915.2 & 1911.4 \\
        r1000.1 & 24 & 24 & 24 & 24 & 11 & 24 & 24 & 20.0 & 20.0 & 20.0 & 20.0 & 20.017 & 20.0 & 20.0 & 210.3 & 263.0 & 145.4 & 224.3 & 1473.1 & 633.1 & 300.2 \\
        DSJC1000.5 & 0 & 0 & 0 & 0 & 0 & 0 & 0 & 78.028 & 78.401 & 90.751 & 90.786 & 88.568 & 84.175 & 87.938 & 1841.3 & 1869.4 & 1825.2 & 1836.2 & 1850.5 & 1837.6 & 1828.3 \\
        flat1000\_60\_0 & 0 & 0 & 0 & 0 & 0 & 0 & 0 & 76.308 & 75.755 & 89.235 & 89.272 & 87.034 & 81.522 & 85.709 & 1841.7 & 1867.8 & 1838.2 & 1828.5 & 1851.6 & 1836.4 & 1821.9 \\
        flat1000\_76\_0 & 0 & 0 & 0 & 0 & 0 & 0 & 0 & 77.293 & 77.134 & 89.58 & 89.589 & 87.434 & 82.189 & 86.803 & 1844.9 & 1861.5 & 1830.5 & 1829.8 & 1835.0 & 1827.9 & 1823.8 \\
        flat1000\_50\_0 & 0 & 0 & 0 & 0 & 0 & 0 & 0 & 51.76 & 61.696 & 88.794 & 88.788 & 86.42 & 66.964 & 80.768 & 1828.8 & 1845.9 & 1830.5 & 1838.6 & 1850.3 & 1828.0 & 1815.7 \\
        5-FullIns\_4 & 10 & 0 & 0 & 18 & 0 & 0 & 0 & 7.283 & 7.616 & 7.407 & 7.285 & 7.715 & 7.325 & 7.902 & 1584.0 & 1827.4 & 1823.4 & 1567.2 & 2084.8 & 1833.0 & 1807.3 \\
        will199GPIA & 23 & 0 & 5 & 5 & 0 & 0 & 0 & 6.2 & 6.244 & 6.204 & 6.208 & 6.543 & 6.241 & 6.482 & 818.2 & 1833.6 & 1787.0 & 1793.6 & 2171.2 & 1829.5 & 2422.4 \\
        wap05a & 24 & 0 & 24 & 24 & 24 & 0 & 0 & 50.0 & 50.396 & 50.0 & 50.0 & 50.0 & 55.367 & 51.529 & 80.2 & 1837.7 & 426.0 & 365.7 & 435.4 & 1829.1 & 1812.5 \\
        wap06a & 12 & 0 & 0 & 0 & 0 & 0 & 0 & 40.012 & 51.262 & 43.622 & 43.47 & 42.989 & 56.945 & 53.552 & 1770.2 & 1825.1 & 1823.2 & 1841.5 & 1868.6 & 1816.9 & 1811.4 \\
        DSJC1000.9 & 24 & 24 & 0 & 0 & 24 & 0 & 0 & 214.855 & 214.855 & 228.839 & 228.876 & 214.855 & 214.973 & 219.215 & 274.1 & 372.9 & 1813.4 & 1831.6 & 1579.4 & 1806.1 & 1805.1 \\
        DSJC500.1 & 0 & 0 & 0 & 0 & 0 & 0 & 0 & 11.415 & 12.141 & 12.8 & 12.93 & 13.394 & 11.363 & 11.848 & 1811.2 & 1819.3 & 1821.5 & 1813.1 & 2020.6 & 1824.9 & 1815.7 \\
        2-FullIns\_5 & 0 & 0 & 0 & 0 & 0 & 0 & 0 & 4.737 & 5.61 & 5.247 & 4.829 & 5.593 & 4.973 & 6.105 & 1825.7 & 1820.2 & 1810.8 & 1832.4 & 1845.6 & 1822.9 & 1803.5 \\
        4-Insertions\_4 & 0 & 0 & 0 & 0 & 0 & 0 & 0 & 2.592 & 2.745 & 2.852 & 2.777 & 3.124 & 2.641 & 3.025 & 1821.2 & 1820.1 & 1809.2 & 1833.7 & 1836.5 & 1832.5 & 1803.5 \\
        2-Insertions\_5 & 0 & 0 & 0 & 0 & 0 & 0 & 0 & 2.92 & 3.458 & 3.349 & 3.22 & 3.73 & 2.978 & 3.747 & 1856.7 & 1823.3 & 1813.4 & 1834.1 & 1828.5 & 1841.8 & 1803.4 \\
        4-FullIns\_4 & 24 & 0 & 24 & 17 & 0 & 22 & 0 & 6.329 & 6.389 & 6.329 & 6.33 & 6.432 & 6.329 & 6.575 & 636.0 & 1817.5 & 778.0 & 1605.9 & 1922.9 & 1099.2 & 1803.2 \\
        r1000.5 & 14 & 24 & 0 & 0 & 24 & 0 & 0 & 234.053 & 234.0 & 268.258 & 268.283 & 234.0 & 249.398 & 243.494 & 1590.5 & 1155.1 & 1820.8 & 1819.0 & 786.8 & 1816.5 & 1803.1 \\
        DSJC500.5 & 0 & 0 & 0 & 0 & 0 & 0 & 0 & 42.453 & 42.628 & 50.549 & 50.607 & 42.347 & 43.559 & 44.098 & 1804.5 & 1815.7 & 1815.8 & 1815.3 & 1808.1 & 1805.3 & 1801.8 \\
        1-Insertions\_6 & 0 & 0 & 0 & 0 & 0 & 0 & 0 & 3.227 & 4.264 & 3.918 & 3.757 & 4.029 & 3.324 & 4.664 & 1830.3 & 1820.3 & 1815.1 & 1821.4 & 1830.1 & 1822.5 & 1802.1 \\
        le450\_5a & 0 & 0 & 0 & 0 & 0 & 0 & 0 & 6.247 & 6.031 & 7.017 & 6.927 & 8.88 & 6.129 & 7.336 & 1813.9 & 1816.9 & 1820.0 & 1820.5 & 1936.5 & 1818.7 & 1801.9 \\
        le450\_5b & 0 & 0 & 0 & 0 & 0 & 0 & 0 & 6.299 & 5.9 & 7.025 & 6.938 & 8.837 & 6.113 & 7.321 & 1810.6 & 1820.3 & 1822.0 & 1824.6 & 1910.7 & 1830.5 & 1802.1 \\
        r1000.1c & 24 & 24 & 0 & 0 & 24 & 24 & 1 & 95.057 & 95.057 & 95.197 & 95.2 & 95.057 & 95.057 & 95.071 & 21.5 & 128.8 & 1828.0 & 1828.6 & 68.7 & 53.2 & 1800.3 \\
        le450\_25a & 24 & 24 & 24 & 24 & 24 & 24 & 24 & 25.0 & 25.0 & 25.0 & 25.0 & 25.0 & 25.0 & 25.0 & 11.4 & 90.8 & 11.5 & 155.0 & 10.5 & 637.3 & 98.3 \\
        le450\_15d & 0 & 0 & 0 & 0 & 0 & 0 & 0 & 17.725 & 21.902 & 20.015 & 19.318 & 21.687 & 21.862 & 21.212 & 1807.4 & 1822.9 & 1826.7 & 1813.2 & 1822.8 & 1812.7 & 1802.7 \\
        le450\_15c & 0 & 0 & 0 & 0 & 0 & 0 & 0 & 17.553 & 21.913 & 20.061 & 19.304 & 21.458 & 21.734 & 21.002 & 1807.3 & 1816.0 & 1821.7 & 1812.8 & 1826.6 & 1811.5 & 1802.8 \\
        le450\_15b & 24 & 0 & 24 & 24 & 0 & 0 & 0 & 15.0 & 17.082 & 15.0 & 15.0 & 15.932 & 17.327 & 15.921 & 109.4 & 1816.1 & 1265.6 & 970.5 & 1895.9 & 1814.9 & 1802.2 \\
        le450\_15a & 24 & 0 & 24 & 24 & 0 & 0 & 0 & 15.0 & 17.186 & 15.0 & 15.0 & 15.87 & 17.557 & 15.941 & 127.6 & 1820.1 & 1319.0 & 977.5 & 1872.4 & 1812.1 & 1801.6 \\
        le450\_25c & 24 & 0 & 0 & 0 & 0 & 0 & 0 & 25.0 & 28.468 & 25.428 & 25.166 & 25.942 & 29.577 & 27.268 & 254.0 & 1818.9 & 1812.3 & 1827.4 & 1827.5 & 1810.2 & 1800.8 \\
        le450\_25d & 24 & 0 & 0 & 0 & 0 & 0 & 0 & 25.0 & 28.011 & 25.408 & 25.136 & 26.091 & 29.055 & 27.246 & 188.5 & 1820.2 & 1807.9 & 1824.6 & 1826.7 & 1812.9 & 1801.2 \\
        le450\_5c & 0 & 7 & 0 & 0 & 0 & 0 & 0 & 5.346 & 5.109 & 7.423 & 7.112 & 7.477 & 5.44 & 7.071 & 1805.3 & 1462.8 & 1819.5 & 1815.1 & 1821.3 & 1817.6 & 1803.4 \\
        \bottomrule
    \end{tabular}}
    \caption{Solving statistics of CG. All methods run on a single CPU core given a cutoff time of $1800$ seconds. For a graph, the first category of columns shows the number of solved test instances out of $24$ instances generated using that graph. The rest two categories respectively show the mean LP objective and mean computational time.}
 \label{tab:cg-solving1}
\end{table}
\end{landscape}

\begin{landscape}
\begin{table}[h!]
    \centering
    \resizebox{\linewidth}{!}{\begin{tabular}{@{}lrrrrrrr|ccccccc|rrrrrrr@{}}
        \toprule
        \multirow{2}{*}{Graph} & \multicolumn{7}{c}{\# solve runs} & \multicolumn{7}{c}{LP objective value} & \multicolumn{7}{c}{Time} \\
         & MLPH & ACO & Gurobi & Gurobi-heur& TSM & Fastwclq & LSCC & MLPH & ACO & Gurobi & Gurobi-heur & TSM & Fastwclq & LSCC & MLPH & ACO & Gurobi & Gurobi-heur & TSM & Fastwclq & LSCC\\ 
        \cmidrule(lr){1-1}\cmidrule(lr){2-8}\cmidrule(lr){9-15}\cmidrule(lr){16-22}
        queen16\_16 & 24 & 24 & 24 & 24 & 0 & 17 & 24 & 16.0 & 16.0 & 16.0 & 16.0 & 16.805 & 16.07 & 16.0 & 34.8 & 135.0 & 418.2 & 506.5 & 1818.5 & 830.4 & 697.8 \\
        le450\_5d & 0 & 11 & 0 & 0 & 0 & 0 & 0 & 5.343 & 5.092 & 7.395 & 6.997 & 7.565 & 5.351 & 7.037 & 1808.1 & 1311.5 & 1816.8 & 1816.6 & 1823.0 & 1814.0 & 1805.4 \\
        DSJC500.9 & 24 & 24 & 0 & 0 & 24 & 24 & 24 & 122.306 & 122.306 & 125.304 & 125.23 & 122.306 & 122.306 & 122.306 & 12.6 & 168.4 & 1824.1 & 1827.0 & 71.7 & 111.6 & 389.0 \\
        3-FullIns\_4 & 24 & 24 & 24 & 24 & 24 & 24 & 23 & 5.392 & 5.392 & 5.392 & 5.392 & 5.392 & 5.392 & 5.392 & 49.4 & 533.3 & 119.7 & 789.1 & 327.7 & 220.2 & 704.6 \\
        3-Insertions\_4 & 0 & 0 & 0 & 0 & 0 & 0 & 0 & 2.533 & 2.619 & 2.604 & 2.715 & 2.79 & 2.477 & 2.796 & 1813.8 & 1818.3 & 1805.0 & 1819.6 & 1819.9 & 1811.3 & 1801.1 \\
        queen15\_15 & 24 & 24 & 24 & 24 & 0 & 22 & 24 & 15.0 & 15.0 & 15.0 & 15.0 & 15.326 & 15.092 & 15.0 & 18.9 & 58.1 & 269.0 & 399.8 & 1824.6 & 319.0 & 373.8 \\
        school1 & 0 & 0 & 0 & 0 & 0 & 0 & 0 & 15.42 & 21.96 & 19.523 & 18.627 & 15.87 & 16.188 & 18.333 & 1804.1 & 1811.1 & 1813.0 & 1819.1 & 1803.1 & 1804.7 & 1800.6 \\
        flat300\_26\_0 & 0 & 0 & 0 & 0 & 0 & 0 & 0 & 26.0 & 26.0 & 30.654 & 31.302 & 26.0 & 26.254 & 26.0 & 1802.6 & 1806.8 & 1828.8 & 1816.9 & 1801.6 & 1800.8 & 1800.4 \\
        flat300\_28\_0 & 0 & 9 & 0 & 0 & 24 & 19 & 3 & 27.525 & 27.52 & 30.78 & 31.266 & 27.52 & 27.698 & 27.52 & 1801.1 & 1678.0 & 1828.0 & 1818.7 & 225.0 & 1203.9 & 1761.4 \\
        DSJR500.5 & 24 & 24 & 0 & 0 & 24 & 21 & 24 & 122.0 & 122.0 & 132.687 & 132.97 & 122.0 & 122.399 & 122.0 & 16.4 & 159.8 & 1810.8 & 1807.6 & 38.9 & 677.9 & 274.4 \\
        school1\_nsh & 0 & 0 & 0 & 0 & 0 & 0 & 0 & 15.657 & 21.559 & 18.142 & 17.589 & 15.176 & 17.25 & 17.854 & 1803.1 & 1810.7 & 1814.7 & 1811.4 & 1802.8 & 1802.8 & 1800.6 \\
        flat300\_20\_0 & 0 & 0 & 0 & 0 & 0 & 0 & 0 & 20.0 & 20.0 & 24.526 & 29.374 & 20.522 & 20.0 & 20.896 & 1802.9 & 1806.1 & 1816.7 & 1824.8 & 1801.5 & 1801.8 & 1800.4 \\
        DSJC250.1 & 0 & 0 & 0 & 0 & 0 & 0 & 0 & 7.079 & 7.625 & 7.466 & 7.499 & 8.033 & 7.004 & 7.33 & 1973.5 & 1817.1 & 1813.0 & 1808.8 & 1821.8 & 1804.5 & 1833.9 \\
        queen14\_14 & 24 & 24 & 24 & 24 & 24 & 23 & 24 & 14.0 & 14.0 & 14.0 & 14.0 & 14.0 & 14.017 & 14.0 & 12.1 & 31.2 & 183.7 & 343.4 & 220.1 & 235.0 & 218.8 \\
        DSJC250.5 & 5 & 24 & 0 & 0 & 24 & 24 & 17 & 25.166 & 25.165 & 26.714 & 27.118 & 25.165 & 25.165 & 25.165 & 1746.4 & 310.4 & 1816.2 & 1807.5 & 83.9 & 302.2 & 1240.2 \\
        queen13\_13 & 24 & 24 & 24 & 24 & 24 & 24 & 24 & 13.0 & 13.0 & 13.0 & 13.0 & 13.0 & 13.0 & 13.0 & 5.5 & 13.4 & 91.3 & 245.5 & 42.9 & 47.0 & 104.8 \\
        DSJR500.1c & 24 & 24 & 11 & 13 & 24 & 24 & 24 & 84.136 & 84.136 & 84.158 & 84.161 & 84.136 & 84.136 & 84.136 & 1.8 & 62.2 & 1152.5 & 1171.6 & 1.2 & 1.4 & 26.5 \\
        1-FullIns\_5 & 24 & 1 & 15 & 0 & 0 & 24 & 0 & 3.909 & 3.957 & 3.91 & 3.933 & 3.939 & 3.909 & 3.974 & 13.0 & 1777.8 & 1581.4 & 1813.2 & 1808.5 & 118.7 & 1800.3 \\
        1-Insertions\_5 & 0 & 0 & 0 & 0 & 0 & 24 & 0 & 2.953 & 3.13 & 2.963 & 3.11 & 2.988 & 2.943 & 3.149 & 1801.5 & 1813.2 & 1803.2 & 1816.7 & 1803.8 & 449.8 & 1800.4 \\
        queen12\_12 & 24 & 24 & 24 & 24 & 24 & 24 & 24 & 12.0 & 12.0 & 12.0 & 12.0 & 12.0 & 12.0 & 12.0 & 3.4 & 7.2 & 60.2 & 217.4 & 13.3 & 29.1 & 60.3 \\
        2-FullIns\_4 & 24 & 24 & 24 & 24 & 24 & 24 & 24 & 4.485 & 4.485 & 4.485 & 4.485 & 4.485 & 4.485 & 4.485 & 4.0 & 157.3 & 14.0 & 391.5 & 15.6 & 17.9 & 298.0 \\
        DSJC250.9 & 24 & 24 & 0 & 0 & 24 & 24 & 24 & 70.392 & 70.392 & 71.24 & 71.242 & 70.392 & 70.392 & 70.392 & 1.5 & 6.2 & 1833.7 & 1836.3 & 4.9 & 6.5 & 16.4 \\
        2-Insertions\_4 & 0 & 0 & 24 & 0 & 24 & 24 & 0 & 2.567 & 2.592 & 2.56 & 2.671 & 2.56 & 2.56 & 2.6 & 1801.0 & 1803.1 & 439.7 & 1805.8 & 1328.3 & 239.0 & 1800.5 \\
        r250.5 & 24 & 24 & 24 & 24 & 24 & 24 & 24 & 65.0 & 65.0 & 65.0 & 65.0 & 65.0 & 65.0 & 65.0 & 2.2 & 11.3 & 476.3 & 282.3 & 4.3 & 61.4 & 26.1 \\
        queen11\_11 & 24 & 24 & 24 & 24 & 24 & 24 & 24 & 11.0 & 11.0 & 11.0 & 11.0 & 11.0 & 11.0 & 11.0 & 2.0 & 3.6 & 36.3 & 182.5 & 5.9 & 14.5 & 32.0 \\
        mug100\_1 & 24 & 24 & 24 & 0 & 24 & 24 & 24 & 3.03 & 3.03 & 3.03 & 3.044 & 3.03 & 3.03 & 3.03 & 9.5 & 27.5 & 17.5 & 1811.1 & 18.6 & 54.8 & 48.3 \\
        mug100\_25 & 24 & 24 & 24 & 2 & 24 & 24 & 24 & 3.03 & 3.03 & 3.03 & 3.036 & 3.03 & 3.03 & 3.03 & 6.6 & 21.6 & 16.4 & 1809.1 & 15.8 & 48.8 & 42.6 \\
        mug88\_1 & 24 & 24 & 24 & 14 & 24 & 24 & 24 & 3.034 & 3.034 & 3.034 & 3.036 & 3.034 & 3.034 & 3.034 & 3.4 & 9.3 & 8.7 & 1640.1 & 6.5 & 25.7 & 18.8 \\
        myciel7 & 24 & 0 & 24 & 0 & 9 & 24 & 0 & 4.095 & 4.16 & 4.095 & 4.212 & 4.113 & 4.095 & 4.097 & 132.8 & 1807.9 & 1097.1 & 1807.6 & 1625.9 & 124.2 & 1800.2 \\
        DSJC125.1 & 0 & 0 & 24 & 0 & 24 & 0 & 0 & 4.458 & 4.565 & 4.454 & 4.47 & 4.454 & 4.484 & 4.461 & 1800.3 & 1801.1 & 495.2 & 1808.7 & 75.2 & 1800.6 & 1800.2 \\
        mug88\_25 & 24 & 24 & 24 & 5 & 24 & 24 & 24 & 3.034 & 3.034 & 3.034 & 3.038 & 3.034 & 3.034 & 3.034 & 3.9 & 9.3 & 9.4 & 1761.2 & 8.1 & 23.1 & 22.6 \\
        queen10\_10 & 24 & 24 & 24 & 24 & 24 & 24 & 24 & 10.0 & 10.0 & 10.0 & 10.0 & 10.0 & 10.0 & 10.0 & 1.2 & 1.7 & 28.1 & 158.2 & 3.7 & 8.4 & 21.1 \\
        DSJC125.5 & 24 & 24 & 24 & 6 & 24 & 24 & 24 & 15.727 & 15.727 & 15.727 & 15.727 & 15.727 & 15.727 & 15.727 & 56.2 & 7.8 & 908.7 & 1607.6 & 4.8 & 14.3 & 53.5 \\
        4-Insertions\_3 & 24 & 24 & 24 & 24 & 24 & 24 & 24 & 2.276 & 2.276 & 2.276 & 2.276 & 2.276 & 2.276 & 2.276 & 362.8 & 299.8 & 16.7 & 882.1 & 8.8 & 25.6 & 61.5 \\
        queen9\_9 & 24 & 24 & 24 & 24 & 24 & 24 & 24 & 9.0 & 9.0 & 9.0 & 9.0 & 9.0 & 9.0 & 9.0 & 0.5 & 0.7 & 9.8 & 132.9 & 1.4 & 3.5 & 7.0 \\
        r125.5 & 24 & 24 & 24 & 24 & 24 & 24 & 24 & 36.0 & 36.0 & 36.0 & 36.0 & 36.0 & 36.0 & 36.0 & 0.3 & 0.7 & 3.5 & 77.8 & 0.3 & 1.0 & 1.2 \\
        DSJC125.9 & 24 & 24 & 24 & 24 & 24 & 24 & 24 & 42.727 & 42.727 & 42.727 & 42.727 & 42.727 & 42.727 & 42.727 & 0.2 & 0.5 & 53.7 & 115.4 & 0.3 & 0.4 & 0.5 \\
        1-FullIns\_4 & 24 & 24 & 24 & 24 & 24 & 24 & 24 & 3.633 & 3.633 & 3.633 & 3.633 & 3.633 & 3.633 & 3.633 & 0.2 & 1.3 & 2.0 & 71.0 & 0.9 & 0.7 & 1.4 \\
        1-Insertions\_4 & 24 & 24 & 24 & 24 & 24 & 24 & 24 & 2.774 & 2.774 & 2.774 & 2.774 & 2.774 & 2.774 & 2.774 & 3.3 & 39.0 & 9.2 & 484.7 & 2.4 & 4.5 & 19.2 \\
        myciel6 & 24 & 24 & 24 & 24 & 24 & 24 & 24 & 3.834 & 3.834 & 3.834 & 3.834 & 3.834 & 3.834 & 3.834 & 1.9 & 518.7 & 28.6 & 794.1 & 12.1 & 9.9 & 76.1 \\
        3-Insertions\_3 & 24 & 24 & 24 & 24 & 24 & 24 & 24 & 2.334 & 2.334 & 2.334 & 2.334 & 2.334 & 2.334 & 2.334 & 6.8 & 13.2 & 3.5 & 300.1 & 1.1 & 4.1 & 11.2 \\
        queen8\_8 & 24 & 24 & 24 & 24 & 24 & 24 & 24 & 8.444 & 8.444 & 8.444 & 8.444 & 8.444 & 8.444 & 8.444 & 0.1 & 0.1 & 1.3 & 60.1 & 0.2 & 0.4 & 0.4 \\
        2-Insertions\_3 & 24 & 24 & 24 & 24 & 24 & 24 & 24 & 2.423 & 2.423 & 2.423 & 2.423 & 2.423 & 2.423 & 2.423 & 0.2 & 0.3 & 0.6 & 112.8 & 0.2 & 0.4 & 0.4 \\
        myciel5 & 24 & 24 & 24 & 24 & 24 & 24 & 24 & 3.553 & 3.553 & 3.553 & 3.553 & 3.553 & 3.553 & 3.553 & 0.1 & 1.2 & 1.1 & 139.0 & 0.3 & 0.7 & 1.6 \\
        myciel4 & 24 & 24 & 24 & 24 & 24 & 24 & 24 & 3.245 & 3.245 & 3.245 & 3.245 & 3.245 & 3.245 & 3.245 & 0.0 & 0.0 & 0.1 & 2.2 & 0.0 & 0.0 & 0.0 \\
        \bottomrule
    \end{tabular}}
     \caption{Solving statistics of CG (Table \ref{tab:cg-solving1} continued). }
 \label{tab:cg-solving2}
\end{table}
\end{landscape}

\begin{landscape}
\begin{table}[phtb]
    \centering
    \resizebox{\linewidth}{!}{\begin{tabular}{@{}lrrrrrrr|lllllll|lllllll@{}}
        \toprule
                \multirow{2}{*}{Graph} & \multicolumn{7}{c}{\# columns with negative reduced costs} & \multicolumn{7}{c}{Minimum reduced cost} & \multicolumn{7}{c}{ Mean reduced cost} \\
         & MLPH & ACO & Gurobi & Gurobi-heur& TSM & Fastwclq & LSCC & MLPH & ACO & Gurobi & Gurobi-heur & TSM & Fastwclq & LSCC & MLPH & ACO & Gurobi & Gurobi-heur & TSM & Fastwclq & LSCC\\ 
        \cmidrule(lr){1-1}\cmidrule(lr){2-8}\cmidrule(lr){9-15}\cmidrule(lr){16-22}
        ash608GPIA & 15254.2 & 277.1 & 26.9 & 438.7 & 4695.6 & 22.3 & 0.9 & -0.11 & -0.05 & -0.18 & -0.18 & -0.11 & -0.06 & -0.02 & -0.06 & -0.01 & -0.13 & -0.08 & -0.07 & -0.02 & -0.02 \\
        wap08a & 57150.2 & 206.5 & 3.8 & 6.5 & 509.2 & 0.2 & 1.0 & -1.59 & -0.37 & -1.69 & -1.91 & -0.56 & -0.01 & -0.55 & -0.99 & -0.07 & -1.17 & -1.13 & -0.27 & -0.01 & -0.55 \\
        wap07a & 63284.2 & 249.2 & 3.2 & 6.1 & 493.2 & 0.7 & 1.0 & -1.55 & -0.37 & -1.47 & -1.83 & -0.59 & -0.02 & -0.75 & -0.98 & -0.07 & -1.0 & -1.1 & -0.28 & -0.02 & -0.75 \\
        abb313GPIA & 26178.0 & 142.2 & 13.2 & 321.3 & 1734.8 & 41.0 & 0.9 & -0.45 & -0.16 & -0.63 & -0.63 & -0.46 & -0.23 & -0.16 & -0.19 & -0.03 & -0.44 & -0.32 & -0.27 & -0.06 & -0.16 \\
        3-FullIns\_5 & 5738.9 & 9.5 & 3.0 & 1164.9 & 0.0 & 502.4 & 0.6 & -2.62 & -0.41 & -2.63 & -2.63 & 0.0 & -2.44 & -0.33 & -1.66 & -0.14 & -2.16 & -2.32 & 0.0 & -0.76 & -0.33 \\
        DSJC1000.1 & 50057.5 & 2319.8 & 2.9 & 2.5 & 107.0 & 440.2 & 1.0 & -0.58 & -0.25 & -0.27 & -0.27 & -0.14 & -0.45 & -0.46 & -0.24 & -0.04 & -0.16 & -0.17 & -0.07 & -0.09 & -0.46 \\
        3-Insertions\_5 & 17041.5 & 22.2 & 2.0 & 5128.4 & 5793.0 & 828.8 & 0.8 & -3.19 & -0.38 & -3.24 & -3.24 & -2.22 & -3.0 & -0.18 & -2.88 & -0.1 & -2.84 & -2.38 & -1.12 & -1.33 & -0.18 \\
        r1000.1 & 50030.5 & 1159.0 & 7.3 & 2892.4 & 1749.4 & 52.2 & 1.0 & -0.5 & -0.17 & -0.7 & -0.7 & -0.51 & -0.22 & -0.49 & -0.25 & -0.03 & -0.36 & -0.35 & -0.26 & -0.05 & -0.49 \\
        DSJC1000.5 & 37906.0 & 26552.0 & 2.0 & 1.7 & 23.5 & 298.2 & 1.0 & -0.49 & -0.41 & -0.17 & -0.16 & -0.4 & -0.43 & -0.45 & -0.1 & -0.06 & -0.13 & -0.13 & -0.21 & -0.07 & -0.45 \\
        flat1000\_60\_0 & 34218.0 & 32003.5 & 1.6 & 1.5 & 21.7 & 412.5 & 1.0 & -0.49 & -0.49 & -0.16 & -0.16 & -0.39 & -0.45 & -0.51 & -0.09 & -0.06 & -0.14 & -0.15 & -0.2 & -0.06 & -0.51 \\
        flat1000\_76\_0 & 35215.0 & 33548.5 & 1.4 & 1.4 & 23.0 & 401.1 & 1.0 & -0.46 & -0.42 & -0.18 & -0.19 & -0.4 & -0.45 & -0.45 & -0.09 & -0.06 & -0.16 & -0.17 & -0.21 & -0.06 & -0.45 \\
        flat1000\_50\_0 & 32684.0 & 19579.7 & 1.6 & 1.5 & 25.6 & 236.0 & 1.0 & -0.72 & -0.87 & -0.22 & -0.21 & -0.44 & -0.65 & -0.59 & -0.09 & -0.06 & -0.17 & -0.18 & -0.22 & -0.07 & -0.59 \\
        5-FullIns\_4 & 17277.5 & 32.0 & 3.6 & 4770.3 & 449.9 & 337.1 & 0.6 & -0.92 & -0.24 & -0.96 & -0.96 & -0.96 & -0.77 & -0.12 & -0.43 & -0.06 & -0.76 & -0.94 & -0.47 & -0.12 & -0.12 \\
        will199GPIA & 35012.4 & 615.0 & 7.3 & 8904.2 & 1290.8 & 151.2 & 1.0 & -0.5 & -0.19 & -0.57 & -0.57 & -0.48 & -0.36 & -0.17 & -0.25 & -0.04 & -0.33 & -0.49 & -0.24 & -0.08 & -0.17 \\
        wap05a & 45304.6 & 291.7 & 8.0 & 214.0 & 401.3 & 2.2 & 1.0 & -1.14 & -0.39 & -1.4 & -1.37 & -1.3 & -0.11 & -0.55 & -0.59 & -0.07 & -1.12 & -0.72 & -0.51 & -0.08 & -0.55 \\
        wap06a & 47401.9 & 231.3 & 7.7 & 207.0 & 385.2 & 1.8 & 1.0 & -1.21 & -0.37 & -1.44 & -1.42 & -1.3 & -0.06 & -0.53 & -0.64 & -0.07 & -1.21 & -0.74 & -0.55 & -0.04 & -0.53 \\
        DSJC1000.9 & 2239.2 & 4207.4 & 1.3 & 1.3 & 11.1 & 12.5 & 1.0 & -0.4 & -0.41 & -0.1 & -0.11 & -0.41 & -0.4 & -0.39 & -0.05 & -0.05 & -0.08 & -0.09 & -0.21 & -0.16 & -0.39 \\
        DSJC500.1 & 25029.4 & 1888.5 & 10.2 & 26.6 & 169.2 & 2463.9 & 1.0 & -0.64 & -0.29 & -0.55 & -0.55 & -0.32 & -0.56 & -0.39 & -0.28 & -0.05 & -0.33 & -0.17 & -0.15 & -0.09 & -0.39 \\
        2-FullIns\_5 & 28020.0 & 58.8 & 2.8 & 5427.5 & 410.6 & 1643.5 & 0.7 & -1.37 & -0.45 & -1.38 & -1.38 & -1.38 & -0.94 & -0.13 & -0.81 & -0.1 & -1.18 & -1.32 & -0.71 & -0.41 & -0.13 \\
        4-Insertions\_4 & 21038.5 & 5408.3 & 2.4 & 17811.8 & 2015.5 & 18527.9 & 1.0 & -0.89 & -0.39 & -0.95 & -0.95 & -0.76 & -0.93 & -0.3 & -0.72 & -0.06 & -0.67 & -0.6 & -0.4 & -0.57 & -0.3 \\
        2-Insertions\_5 & 22637.8 & 187.1 & 2.2 & 11643.6 & 2397.5 & 14467.7 & 1.0 & -1.84 & -0.46 & -1.85 & -1.85 & -1.85 & -1.85 & -0.25 & -1.4 & -0.07 & -1.27 & -1.28 & -0.89 & -1.13 & -0.25 \\
        4-FullIns\_4 & 27861.5 & 244.4 & 3.2 & 11288.0 & 272.7 & 377.7 & 0.6 & -0.6 & -0.27 & -0.63 & -0.63 & -0.63 & -0.57 & -0.14 & -0.27 & -0.06 & -0.46 & -0.61 & -0.31 & -0.1 & -0.14 \\
        r1000.5 & 16634.6 & 13192.6 & 1.9 & 1.7 & 34.8 & 116.5 & 1.0 & -0.6 & -0.5 & -0.25 & -0.26 & -0.67 & -0.45 & -0.55 & -0.12 & -0.08 & -0.2 & -0.21 & -0.34 & -0.08 & -0.55 \\
        DSJC500.5 & 12686.2 & 10771.2 & 2.0 & 2.1 & 23.3 & 302.7 & 1.0 & -0.44 & -0.4 & -0.22 & -0.21 & -0.46 & -0.45 & -0.42 & -0.1 & -0.06 & -0.19 & -0.18 & -0.24 & -0.07 & -0.42 \\
        1-Insertions\_6 & 16322.5 & 73.6 & 2.4 & 7568.8 & 1239.3 & 10858.0 & 0.7 & -2.36 & -0.48 & -2.36 & -2.36 & -2.36 & -2.36 & -0.13 & -1.51 & -0.08 & -1.55 & -1.69 & -1.11 & -0.94 & -0.13 \\
        le450\_5a & 22526.3 & 12282.5 & 16.7 & 2089.0 & 474.0 & 2643.7 & 1.0 & -1.01 & -0.7 & -1.23 & -1.22 & -0.62 & -1.03 & -0.6 & -0.47 & -0.22 & -0.72 & -0.47 & -0.3 & -0.15 & -0.6 \\
        le450\_5b & 22527.8 & 12263.1 & 16.4 & 1630.8 & 481.0 & 2720.1 & 1.0 & -0.97 & -0.64 & -1.2 & -1.2 & -0.58 & -1.07 & -0.54 & -0.43 & -0.2 & -0.64 & -0.4 & -0.28 & -0.17 & -0.54 \\
        r1000.1c & 23.9 & 24.0 & 0.3 & 0.2 & 9.6 & 4.0 & 0.2 & -0.3 & -0.3 & -0.05 & -0.02 & -0.3 & -0.3 & -0.04 & -0.09 & -0.09 & -0.05 & -0.02 & -0.14 & -0.2 & -0.04 \\
        le450\_25a & 22515.0 & 405.4 & 6.9 & 8886.5 & 209.6 & 14.7 & 1.0 & -0.67 & -0.27 & -0.72 & -0.72 & -0.72 & -0.28 & -0.28 & -0.33 & -0.05 & -0.46 & -0.69 & -0.3 & -0.1 & -0.28 \\
        le450\_15d & 22474.0 & 766.8 & 10.8 & 366.9 & 180.8 & 105.2 & 1.0 & -1.09 & -0.43 & -1.13 & -1.11 & -1.17 & -0.61 & -0.47 & -0.55 & -0.07 & -0.71 & -0.33 & -0.49 & -0.11 & -0.47 \\
        le450\_15c & 22496.5 & 740.5 & 11.1 & 385.7 & 164.7 & 126.9 & 1.0 & -1.09 & -0.42 & -1.12 & -1.13 & -1.17 & -0.63 & -0.45 & -0.55 & -0.07 & -0.71 & -0.35 & -0.49 & -0.11 & -0.45 \\
        le450\_15b & 22523.9 & 227.6 & 11.2 & 3345.6 & 291.1 & 49.7 & 1.0 & -1.03 & -0.32 & -1.16 & -1.16 & -1.16 & -0.46 & -0.54 & -0.57 & -0.07 & -0.8 & -0.5 & -0.47 & -0.1 & -0.54 \\
        le450\_15a & 22523.9 & 247.5 & 10.8 & 3329.4 & 290.8 & 68.3 & 1.0 & -1.07 & -0.34 & -1.18 & -1.18 & -1.18 & -0.5 & -0.51 & -0.58 & -0.07 & -0.81 & -0.51 & -0.47 & -0.1 & -0.51 \\
        le450\_25c & 22429.9 & 541.1 & 10.0 & 767.9 & 129.0 & 28.0 & 1.0 & -1.02 & -0.41 & -1.06 & -1.05 & -1.06 & -0.46 & -0.42 & -0.52 & -0.07 & -0.67 & -0.35 & -0.44 & -0.13 & -0.42 \\
        le450\_25d & 22423.7 & 690.7 & 9.7 & 882.8 & 133.6 & 32.2 & 1.0 & -0.98 & -0.39 & -1.03 & -1.01 & -1.03 & -0.47 & -0.34 & -0.49 & -0.07 & -0.69 & -0.37 & -0.43 & -0.12 & -0.34 \\
        le450\_5c & 22504.5 & 27218.2 & 10.5 & 4451.0 & 420.9 & 4385.6 & 1.0 & -1.33 & -1.25 & -1.43 & -1.43 & -1.43 & -1.42 & -0.73 & -0.48 & -0.44 & -0.73 & -1.07 & -1.12 & -0.35 & -0.73 \\
        \bottomrule
    \end{tabular}}
     \caption{Statistics of solving the MWISP at the initial iteration of CG. The first category shows the number of found columns with negative reduced costs. The second category shows the most negative reduced cost of the best-found solution for MWISP. The third category shows the mean reduced cost among the columns with negative reduced costs. For a graph, the statistics are averaged over $24$ problem instances generated using that graph.}
 \label{tab:cg-pp1}
\end{table}
\end{landscape}

\begin{landscape}
\begin{table}[phtb]
    \centering
    \resizebox{\linewidth}{!}{\begin{tabular}{@{}lrrrrrrr|lllllll|lllllll@{}}
        \toprule
                \multirow{2}{*}{Graph} & \multicolumn{7}{c}{\# columns with negative reduced costs} & \multicolumn{7}{c}{Minimum reduced cost} & \multicolumn{7}{c}{ Mean reduced cost} \\
         & MLPH & ACO & Gurobi & Gurobi-heur& TSM & Fastwclq & LSCC & MLPH & ACO & Gurobi & Gurobi-heur & TSM & Fastwclq & LSCC & MLPH & ACO & Gurobi & Gurobi-heur & TSM & Fastwclq & LSCC\\ 
        \cmidrule(lr){1-1}\cmidrule(lr){2-8}\cmidrule(lr){9-15}\cmidrule(lr){16-22}
         
        queen16\_16 & 7207.5 & 13696.4 & 9.9 & 1588.5 & 27.6 & 112.9 & 1.0 & -0.23 & -0.26 & -0.3 & -0.3 & -0.27 & -0.24 & -0.26 & -0.05 & -0.05 & -0.21 & -0.1 & -0.14 & -0.04 & -0.26 \\
        le450\_5d & 22511.0 & 28603.0 & 13.0 & 4793.6 & 851.1 & 4660.2 & 1.0 & -1.3 & -1.27 & -1.39 & -1.39 & -1.39 & -1.37 & -0.75 & -0.47 & -0.42 & -0.77 & -1.11 & -0.72 & -0.32 & -0.75 \\
        DSJC500.9 & 437.2 & 489.9 & 0.5 & 0.5 & 7.7 & 8.3 & 1.0 & -0.41 & -0.41 & -0.04 & -0.06 & -0.39 & -0.41 & -0.41 & -0.07 & -0.07 & -0.04 & -0.06 & -0.21 & -0.19 & -0.41 \\
        3-FullIns\_4 & 13700.9 & 1340.0 & 3.8 & 21490.0 & 124.0 & 264.1 & 0.8 & -0.44 & -0.3 & -0.45 & -0.45 & -0.45 & -0.42 & -0.13 & -0.17 & -0.07 & -0.35 & -0.45 & -0.23 & -0.08 & -0.13 \\
        3-Insertions\_4 & 8786.8 & 16248.0 & 2.6 & 24684.9 & 1159.6 & 7718.2 & 1.0 & -0.73 & -0.39 & -0.77 & -0.77 & -0.77 & -0.75 & -0.31 & -0.49 & -0.06 & -0.52 & -0.55 & -0.38 & -0.36 & -0.31 \\
        queen15\_15 & 5518.4 & 10581.1 & 9.9 & 1706.4 & 28.8 & 93.3 & 1.0 & -0.22 & -0.25 & -0.28 & -0.28 & -0.28 & -0.22 & -0.25 & -0.05 & -0.05 & -0.18 & -0.09 & -0.15 & -0.04 & -0.25 \\
        school1 & 18433.2 & 1238.4 & 5.9 & 470.4 & 140.0 & 47.5 & 1.0 & -1.6 & -1.12 & -1.76 & -1.76 & -1.76 & -1.49 & -0.58 & -0.44 & -0.13 & -1.11 & -0.88 & -0.85 & -0.28 & -0.58 \\
        flat300\_26\_0 & 4947.5 & 4167.2 & 3.6 & 7.1 & 19.9 & 139.9 & 1.0 & -0.44 & -0.41 & -0.23 & -0.18 & -0.45 & -0.43 & -0.4 & -0.08 & -0.06 & -0.16 & -0.09 & -0.22 & -0.08 & -0.4 \\
        flat300\_28\_0 & 5200.3 & 5116.8 & 4.5 & 4.3 & 18.6 & 158.6 & 1.0 & -0.42 & -0.39 & -0.23 & -0.18 & -0.43 & -0.41 & -0.39 & -0.08 & -0.06 & -0.14 & -0.09 & -0.22 & -0.07 & -0.39 \\
        DSJR500.5 & 7059.4 & 6738.0 & 1.7 & 1.5 & 27.8 & 242.0 & 1.0 & -0.48 & -0.48 & -0.14 & -0.14 & -0.57 & -0.44 & -0.48 & -0.09 & -0.08 & -0.11 & -0.12 & -0.31 & -0.06 & -0.48 \\
        school1\_nsh & 16043.3 & 2098.6 & 9.3 & 204.8 & 127.7 & 66.5 & 1.0 & -1.31 & -0.82 & -1.44 & -1.42 & -1.44 & -1.0 & -0.49 & -0.36 & -0.11 & -0.8 & -0.44 & -0.7 & -0.19 & -0.49 \\
        flat300\_20\_0 & 5701.5 & 2828.6 & 5.3 & 5.0 & 28.1 & 73.4 & 1.0 & -0.84 & -0.77 & -0.42 & -0.26 & -0.84 & -0.84 & -0.67 & -0.1 & -0.08 & -0.25 & -0.13 & -0.4 & -0.2 & -0.67 \\
        DSJC250.1 & 12482.8 & 1383.2 & 17.0 & 204.7 & 219.0 & 682.3 & 1.0 & -0.65 & -0.32 & -0.68 & -0.69 & -0.7 & -0.5 & -0.35 & -0.28 & -0.06 & -0.42 & -0.15 & -0.31 & -0.09 & -0.35 \\
        queen14\_14 & 4688.3 & 7213.9 & 9.0 & 1615.8 & 22.7 & 61.5 & 1.0 & -0.22 & -0.25 & -0.27 & -0.27 & -0.27 & -0.22 & -0.25 & -0.05 & -0.05 & -0.18 & -0.09 & -0.15 & -0.05 & -0.25 \\
        DSJC250.5 & 2162.5 & 3218.1 & 5.8 & 5.6 & 18.5 & 81.2 & 1.0 & -0.39 & -0.39 & -0.31 & -0.24 & -0.42 & -0.4 & -0.39 & -0.09 & -0.06 & -0.17 & -0.12 & -0.21 & -0.08 & -0.39 \\
        queen13\_13 & 3767.9 & 5989.8 & 8.8 & 1773.2 & 21.2 & 45.1 & 1.0 & -0.23 & -0.24 & -0.27 & -0.27 & -0.27 & -0.24 & -0.24 & -0.05 & -0.05 & -0.18 & -0.08 & -0.14 & -0.06 & -0.24 \\
        DSJR500.1c & 2.1 & 2.1 & 0.0 & 0.0 & 1.5 & 1.1 & 0.2 & -0.42 & -0.42 & -0.04 & 0.0 & -0.42 & -0.42 & -0.07 & -0.29 & -0.29 & -0.04 & 0.0 & -0.32 & -0.39 & -0.07 \\
        1-FullIns\_5 & 2937.1 & 462.6 & 3.5 & 22120.5 & 132.6 & 1784.4 & 0.9 & -1.02 & -0.58 & -1.02 & -1.02 & -1.02 & -1.02 & -0.68 & -0.46 & -0.14 & -0.71 & -0.98 & -0.54 & -0.27 & -0.68 \\
        1-Insertions\_5 & 3485.7 & 902.0 & 4.8 & 30114.7 & 388.1 & 1221.4 & 0.9 & -0.88 & -0.47 & -0.88 & -0.88 & -0.88 & -0.87 & -0.27 & -0.37 & -0.06 & -0.62 & -0.66 & -0.38 & -0.29 & -0.27 \\
        queen12\_12 & 3135.4 & 3682.1 & 7.4 & 1492.7 & 17.7 & 34.5 & 1.0 & -0.23 & -0.23 & -0.25 & -0.25 & -0.25 & -0.22 & -0.22 & -0.05 & -0.05 & -0.16 & -0.07 & -0.13 & -0.05 & -0.22 \\
        2-FullIns\_4 & 5304.7 & 864.5 & 3.2 & 34333.8 & 102.5 & 92.2 & 0.9 & -0.4 & -0.31 & -0.42 & -0.42 & -0.42 & -0.41 & -0.19 & -0.15 & -0.06 & -0.37 & -0.42 & -0.22 & -0.09 & -0.19 \\
        DSJC250.9 & 83.5 & 84.4 & 0.5 & 0.5 & 4.6 & 4.4 & 1.0 & -0.29 & -0.29 & -0.05 & -0.05 & -0.29 & -0.29 & -0.29 & -0.06 & -0.06 & -0.05 & -0.05 & -0.16 & -0.16 & -0.29 \\
        2-Insertions\_4 & 2593.7 & 3703.5 & 4.0 & 38473.6 & 332.6 & 1536.9 & 1.0 & -0.53 & -0.36 & -0.53 & -0.53 & -0.53 & -0.51 & -0.17 & -0.21 & -0.05 & -0.36 & -0.36 & -0.24 & -0.19 & -0.17 \\
        r250.5 & 2692.6 & 3215.2 & 4.6 & 362.4 & 19.4 & 62.6 & 1.0 & -0.47 & -0.49 & -0.54 & -0.54 & -0.54 & -0.4 & -0.4 & -0.12 & -0.12 & -0.31 & -0.21 & -0.29 & -0.11 & -0.4 \\
        queen11\_11 & 2183.0 & 2510.9 & 7.0 & 1692.8 & 14.6 & 22.4 & 1.0 & -0.22 & -0.22 & -0.24 & -0.24 & -0.24 & -0.23 & -0.21 & -0.05 & -0.04 & -0.15 & -0.07 & -0.12 & -0.07 & -0.21 \\
        mug100\_1 & 4005.9 & 5590.3 & 6.4 & 22685.2 & 186.8 & 796.4 & 1.0 & -0.13 & -0.12 & -0.16 & -0.16 & -0.16 & -0.12 & -0.12 & -0.04 & -0.04 & -0.1 & -0.08 & -0.09 & -0.02 & -0.12 \\
        mug100\_25 & 4399.0 & 5410.9 & 6.1 & 23163.2 & 195.2 & 748.5 & 1.0 & -0.14 & -0.13 & -0.16 & -0.16 & -0.16 & -0.13 & -0.12 & -0.05 & -0.04 & -0.1 & -0.09 & -0.09 & -0.02 & -0.12 \\
        mug88\_1 & 3220.0 & 3931.6 & 6.9 & 19528.0 & 146.7 & 117.4 & 1.0 & -0.11 & -0.11 & -0.13 & -0.13 & -0.13 & -0.11 & -0.1 & -0.03 & -0.03 & -0.09 & -0.07 & -0.07 & -0.02 & -0.1 \\
        myciel7 & 820.5 & 192.4 & 3.4 & 32820.8 & 133.6 & 624.0 & 0.7 & -1.53 & -0.89 & -1.53 & -1.53 & -1.53 & -1.47 & -0.27 & -0.45 & -0.13 & -1.1 & -0.96 & -0.48 & -0.37 & -0.27 \\
        DSJC125.1 & 4929.4 & 1644.3 & 11.0 & 35475.4 & 131.7 & 156.2 & 1.0 & -0.5 & -0.32 & -0.51 & -0.51 & -0.51 & -0.43 & -0.23 & -0.2 & -0.07 & -0.29 & -0.41 & -0.23 & -0.09 & -0.23 \\
        mug88\_25 & 3290.2 & 3953.1 & 5.8 & 22046.8 & 144.0 & 431.0 & 1.0 & -0.12 & -0.12 & -0.15 & -0.15 & -0.15 & -0.12 & -0.12 & -0.04 & -0.03 & -0.1 & -0.07 & -0.08 & -0.02 & -0.12 \\
        queen10\_10 & 1194.7 & 1402.9 & 6.0 & 1502.6 & 10.2 & 14.7 & 1.0 & -0.21 & -0.21 & -0.22 & -0.22 & -0.22 & -0.22 & -0.21 & -0.05 & -0.04 & -0.15 & -0.06 & -0.12 & -0.09 & -0.21 \\
        DSJC125.5 & 325.4 & 620.0 & 7.6 & 123.5 & 11.5 & 14.9 & 1.0 & -0.29 & -0.3 & -0.31 & -0.3 & -0.31 & -0.21 & -0.25 & -0.07 & -0.06 & -0.17 & -0.07 & -0.16 & -0.08 & -0.25 \\
        4-Insertions\_3 & 924.6 & 1832.0 & 2.8 & 33978.0 & 87.2 & 425.2 & 0.9 & -0.23 & -0.22 & -0.24 & -0.24 & -0.24 & -0.24 & -0.09 & -0.07 & -0.05 & -0.16 & -0.1 & -0.13 & -0.06 & -0.09 \\
        queen9\_9 & 545.8 & 610.2 & 4.8 & 1744.5 & 9.7 & 8.9 & 1.0 & -0.21 & -0.21 & -0.22 & -0.22 & -0.22 & -0.18 & -0.2 & -0.05 & -0.05 & -0.14 & -0.04 & -0.12 & -0.08 & -0.2 \\
        r125.5 & 425.1 & 757.3 & 3.3 & 3655.7 & 11.5 & 11.4 & 1.0 & -0.39 & -0.4 & -0.4 & -0.4 & -0.4 & -0.4 & -0.31 & -0.14 & -0.15 & -0.29 & -0.13 & -0.24 & -0.2 & -0.31 \\
        DSJC125.9 & 6.8 & 6.8 & 2.3 & 4.7 & 2.8 & 2.0 & 1.0 & -0.19 & -0.19 & -0.19 & -0.19 & -0.19 & -0.19 & -0.19 & -0.08 & -0.08 & -0.14 & -0.1 & -0.12 & -0.14 & -0.19 \\
        1-FullIns\_4 & 213.7 & 110.3 & 2.9 & 41261.2 & 21.1 & 61.1 & 0.8 & -0.29 & -0.26 & -0.29 & -0.29 & -0.29 & -0.29 & -0.17 & -0.12 & -0.08 & -0.19 & -0.21 & -0.16 & -0.11 & -0.17 \\
        1-Insertions\_4 & 180.7 & 297.4 & 5.3 & 23800.0 & 50.8 & 65.0 & 0.9 & -0.23 & -0.21 & -0.23 & -0.23 & -0.23 & -0.21 & -0.08 & -0.06 & -0.04 & -0.16 & -0.06 & -0.09 & -0.06 & -0.08 \\
        myciel6 & 207.7 & 205.0 & 4.0 & 48811.9 & 68.1 & 134.1 & 0.6 & -0.69 & -0.56 & -0.69 & -0.69 & -0.69 & -0.69 & -0.09 & -0.19 & -0.1 & -0.54 & -0.28 & -0.24 & -0.2 & -0.09 \\
        3-Insertions\_3 & 204.4 & 361.5 & 4.0 & 16366.9 & 34.5 & 100.5 & 1.0 & -0.17 & -0.17 & -0.18 & -0.18 & -0.18 & -0.18 & -0.12 & -0.04 & -0.04 & -0.11 & -0.05 & -0.09 & -0.05 & -0.12 \\
        queen8\_8 & 180.3 & 177.8 & 4.1 & 353.2 & 5.0 & 6.9 & 1.0 & -0.17 & -0.17 & -0.17 & -0.17 & -0.17 & -0.17 & -0.16 & -0.04 & -0.04 & -0.12 & -0.03 & -0.1 & -0.09 & -0.16 \\
        2-Insertions\_3 & 26.3 & 29.5 & 3.1 & 177.9 & 5.8 & 5.9 & 1.0 & -0.07 & -0.07 & -0.07 & -0.07 & -0.07 & -0.07 & -0.07 & -0.03 & -0.03 & -0.06 & -0.02 & -0.04 & -0.04 & -0.07 \\
        myciel5 & 34.0 & 41.2 & 3.5 & 1038.7 & 14.6 & 12.8 & 0.8 & -0.28 & -0.27 & -0.28 & -0.28 & -0.28 & -0.28 & -0.13 & -0.1 & -0.08 & -0.21 & -0.04 & -0.11 & -0.09 & -0.13 \\
        myciel4 & 2.6 & 2.7 & 1.2 & 3.3 & 1.4 & 1.2 & 0.8 & -0.09 & -0.09 & -0.09 & -0.09 & -0.09 & -0.09 & -0.09 & -0.06 & -0.06 & -0.08 & -0.06 & -0.07 & -0.08 & -0.09 \\
        \bottomrule
    \end{tabular}}
     \caption{Statistics of solving the initial MWISP (Table \ref{tab:cg-pp1} continued).}
  \label{tab:cg-pp2}
\end{table}
\end{landscape}

\bibliography{lib}